\documentclass[reqno]{amsart}
\usepackage{amsmath,amsfonts,amssymb,amsthm,enumitem}
\usepackage[usenames,dvipsnames]{xcolor}
\usepackage[colorlinks=true]{hyperref} 
\usepackage{hyperref}
\hypersetup{linkcolor=Violet,citecolor=PineGreen}
\newtheorem{teor}{Theorem}[section]
\newtheorem{lema}[teor]{Lemma}
\newtheorem{prop}[teor]{Proposition}
\newtheorem{coro}[teor]{Corollary}
\theoremstyle{definition}
\newtheorem{defi}[teor]{Definition}

\newtheorem{nota}[teor]{Remark}

\numberwithin{equation}{section}

\newcommand\JM{Mierczy\'nski}

\newcommand{\N}{\mathbb{N}}
\newcommand{\R}{\mathbb{R}}

\newcommand\RR{\ensuremath{\mathbb{R}}}

\newcommand\NN{\ensuremath{\mathbb{N}}}

\newcommand{\PP}{\mathbb{P}}
\newcommand{\be}{\mathbf e}
\newcommand{\w}{\omega}

\newcommand{\nbd}{\nobreakdash}
\newcommand{\n}[1]{\|#1\|}

\newcommand{\norm}[1]{\ensuremath{\lVert#1\rVert}}
\newcommand{\normL}[1]{\ensuremath{\lVert#1\rVert_{L}}}
\newcommand{\normC}[1]{\ensuremath{\lVert#1\rVert_{C}}}
\newcommand{\normAC}[1]{\ensuremath{\lVert#1\rVert_{AC}}}
\newcommand{\lsm}{\left[\begin{smallmatrix}}
\newcommand{\rsm}{\end{smallmatrix}\right]}
\newcommand{\mlsps}{measurable linear skew\nobreakdash-\hspace{0pt}product
semidynamical system}
\DeclareMathOperator{\spanned}{span}
\newcommand{\OFP}{\ensuremath{(\Omega,\mathfrak{F},\PP)}}

\DeclareMathOperator{\lnplus}{ln^{+}}

\DeclareMathOperator{\dist}{dist}
\DeclareMathOperator{\vol}{vol}

\DeclareMathOperator*{\esssup}{ess\,sup}
\begin{document}
\title[Two dynamical approaches to the notion of exponential separation]
{Two dynamical approaches to the notion of exponential separation for random systems of delay differential equations. }
\author[M.~Kryspin]{Marek Kryspin}
\author[J.~Mierczy\'nski]{Janusz Mierczy\'nski}
\address[Marek Kryspin, Janusz Mierczy\'nski]{Faculty of Pure and Applied Mathematics,
Wroc{\l}aw University of Science and Technology,
Wybrze\.ze Wyspia\'nskiego 27, PL-50-370 Wroc{\l}aw, Poland.}
\email[Marek Kryspin]{Marek.Kryspin@pwr.edu.pl}
\email[Janusz Mierczy\'nski]{Janusz.Mierczynski@pwr.edu.pl}
\address[Sylvia Novo, Rafael Obaya]
{Departamento de Matem\'{a}tica Aplicada, Universidad de
Valladolid, Paseo Prado de la Magdalena 3-5, 47011 Valladolid, Spain.}
\email[Sylvia Novo]{sylvia.novo@uva.es}
\email[Rafael Obaya]{rafael.obaya@uva.es}
\thanks{The first author is supported by the National Science Centre, Poland (NCN) under the grant Sonata Bis with a number NCN 2020/38/E/ST1/00153, the second author is partially supported by the Faculty of Pure and Applied Mathematics, Wrocław University of Science and Technology, and the last two authors were partly supported by MICIIN/FEDER project PID2021-125446NB-I00 and by Universidad de Valladolid under project PIP-TCESC-2020}
\author[S.~Novo]{Sylvia Novo}
\author[R.~Obaya]{Rafael Obaya}
\subjclass[2020]{Primary: 37H15, 37L55, 34K06, Secondary: 37A30, 37A40, 37C65, 60H25}
\date{}
\begin{abstract}
 This paper deals with the exponential separation of type II, an important concept for random systems of differential equations with delay, introduced in \JM\ et al.~\cite{MiNoOb1}. Two different approaches to its  existence are presented. The state space $X$ will be a separable ordered Banach space with $\dim X\geq 2$, dual space $X^{*}$ and  positive cone $X^+$ normal and reproducing. In both cases, appropriate cooperativity and irreducibility conditions  are assumed to provide a family of generalized Floquet subspaces. If in addition $X^*$ is also separable, one obtains a exponential separation of type II. When this is not the case, but there is an Oseledets decomposition for the continuous semiflow, the same result holds.
 Detailed examples are given for all the situations, including also a case where the cone is not normal.
\end{abstract}
\keywords{Random dynamical systems, focusing property, generalized principal Floquet bundle, generalized exponential separation, random delay differential systems, Oseledets decomposition}
\maketitle
\section{Introduction}\label{sec-intro}
This paper deals with the existence of principal Floquet subspaces and generalized exponential separation of type II for positive random dynamical systems generated by linear differential equations with finite delay. In particular,  some previous results, obtained in \JM\ et al.~\cite{MiNoOb1},  are now completed and extended to random linear systems of delay differential equations. This reference introduced a new  focusing condition motivated by the non-injectivity of the flow map defined by  a linear delay differential equation, which is also of great interest in the vector case now considered.  In addition to the construction method of the generalized exponential separation from the principal Floquet subspace given in~\cite{MiNoOb1} and now extended to the high dimensional  case, the present paper provides an alternative approach based in the  Oseledets decomposition for the continuous  time flow generated by the system of linear delay  differential equations. This theorem was stated in \JM\ et al.~\cite{MiNoOb2} as an adaptation of the semi-invertible Oseledets theorem  proved in Gonz\'alez-Tolkman and Quas~\cite{GTQu} and Lee~\cite{Lee}. 
\par\smallskip
The top finite Lyapunov exponent of a positive deterministic or random dynamical system is called the \emph{principal Lyapunov exponent} if the associated invariant family of subspaces, where this Lyapunov exponent is reached, is one-dimensional and spanned by a positive vector. In this case the invariant subspace is called the \emph{principal Floquet subspace}.  The \emph{exponential separation} theory was initiated for positive discrete-time deterministic dynamical systems by Ruelle~\cite{Rue0} and  Pol\'{a}\v{c}ik and Tere\v{s}\v{c}\'{a}k~\cite{pote,pote1}.  Subsequently, important contributions for discrete and continuous time flows were obtained by H\'uska and Poláčik~\cite{HusPol},  H\'uska~\cite{Hus}, H\'uska~et al.~\cite{HusPolSaf}, Novo et al.~\cite{NOS2},  \JM\ and Shen~\cite{MiSh1,MiSh3}, and Shen and Yi~\cite{shyi}, among others. In particular,~\cite{NOS2} introduced the exponential separation of type~II, a version of the classical notion adapted to the study of deterministic delay differential equations. The importance of this concept for the study of linear and nonlinear nonautonomous functional differential equations with finite delay can be found in Novo et al.~\cite{noos5}, Calzada et al.~\cite{COS} and Obaya and Sanz~\cite{OS1,OS2}.
\par\smallskip
In the context of random dynamical systems,  $\OFP$  will denote a probability space with an ergodic measure $\PP$, $\theta\colon\R \times \Omega \to\Omega$, $(t,\w) \mapsto \theta_t\w$ is a metric dynamical system, $X$ is an ordered Banach space and $\Phi\colon \R^+ \times \Omega \times X \to \Omega\times X$, $(t, \w, u) \mapsto (\theta_t\w,U_\w (t)\,u)$ a measurable linear skew-product semiflow. The concept of generalized exponential separation  (of type~I) refers to a measurable decomposition $X=E_1(\w) \oplus F_1(\w)$ for $\PP$-a.e.\ $\w\in\Omega$, where $E_1(\w)= \spanned \{ w(\w) \}$ is the principal Floquet subspace, $F_1(\w)$ is an invariant co-dimensional one closed vector subspace that does not admit any strictly positive vector, and  $U_\w(t)$  exhibits an exponential separation on the sum. 
\par\smallskip
Arnold et al.~\cite{Arn}  proved the existence of generalized exponential separation for discrete-time positive random dynamical systems generated by random families of positive matrices. Later, \JM\ and Shen~\cite{MiShPart1} provided the assumptions sufficient for a general random positive linear skew-product semiflows in order to admit a generalized principal Floquet subspace and generalized exponential separation. A wide range of applications of this theory, including random linear skew-product semiflow generated by cooperative families of ordinary differential equations and parabolic partial differential equations, can be found in \JM\ and Shen~\cite{MiShPart2,MiShPart3}. More recently, \JM\ et al.~\cite{MiNoOb1} adapted the previous abstract theory to the case of a non-injective flow map. More precisely, assuming integrability, positivity and a new focusing condition in the terms of the one in~\cite{MiShPart1}, a generalized exponential separation of type II is introduced.  The main difference concerning this new focusing condition is that  it implies the existence of a  positive $T>0$ such that for each $u \in X^+$ and $\w \in \Omega$ then $U_\w (T)\,u=0$ or $U_\w(t)\,u$ is strictly positive. As a consequence of this dichotomy behavior, concerning the measurable decomposition,  now  $F_1(\w) \cap X^+=\{u \in X^+ \colon U_\w(T)\,u=0 \}$, that is,  $F_1(\w)$ can contain positive vectors. The present paper studies  the applicability of this theory to the case of random systems of delay differential equations.
\par\smallskip
The paper is organized into five sections and $X$ will denote a separable ordered Banach space with $\dim X\geq 2$, dual space $X^{*}$ and  positive cone  $X^+$ always reproducing and normal in most of the paper. Section~\ref{sec-prel} collects the main notions, assumptions and results of~\cite{MiNoOb1} used throughout the paper. Assuming integrability, positivity and a focusing condition, Theorem~\ref{Theorem3.10}  asserts that the semiflow $\Phi$ admits a generalized Floquet subspace $E_1(\w)=\spanned\{w(\w)\}$. Moreover, Theorem~\ref{Theorem3.10}(vii) provides the exponential convergence of the normalized trajectories of positive vectors to $w(\w)$ in the  forward and also in the pullback sense. In addition, when $X^*$ is also separable, Theorem~\ref{thm:exponential-separation} concludes that $\Phi$ admits a generalized exponential separation of type II  with measurable decomposition $X=E_1(\w) \oplus F_1(\w)$. The previous properties of exponential convergence assure that the family of subspaces $E_1(\w)$, $ F_1(\w)$ and their associated projections,  can be calculated numerically in applications. This shows the interest of the applications of this part of the theory in the paper.
\par\smallskip
 Section~\ref{sect:Oseledets} assumes that $\Phi$ admits an Oseledets decomposition. As proved in~\cite[Theorem 3.4]{MiNoOb2}  this decomposition exists when   $\Omega$ is  a Lebesgue space and there is a $t_0>0$ such that $U_{\w}(t_0)$ is a compact operator for all $\w \in \Omega$. In particular, the separability of $X^*$ is not required.  If the top Lyapunov exponent $\lambda_{\mathrm{top}} >-\infty$, then  $X=E(\w) \oplus F(\w)$ where $E(\w)=\{ u\in X : \lim_{t\to \infty}\ln \n{U_{\w}(t)\,u}/t = \lambda_{\mathrm{top}} \}$ and 
$F(\w)=\{u\in X : \lim_{t \to \infty}\ln \n{U_\w(t)\,u}/t < \lambda_{\mathrm{top}} \}$ 
for $\PP$-a.e.\ $\w\in\Omega$. Assuming that  $\Phi$ admits a generalized Floquet subspace $E_1(\w)=\spanned\{w(\w)\}$, Theorem~\ref{thm:Appendix} proves that $E_1(\w)=E(\w)$ for $\PP$-a.e. $\w\in\Omega$ and then, the previous Oseledets decomposition provides a generalized exponential separation of type II.
\par\smallskip
Section~\ref{sect:semiflows} deals with systems of  linear random delay differential equations of the form $z'(t)= A(\theta_t \w)\,z(t) + B(\theta_t \w)\,z(t-1)$ where $A$ belongs to $L_1(\Omega)$ and $B$ satisfies  a $L_q$-integrability  condition for some  $1<q<\infty$. Under these assumptions, as shown in \cite{MiNoOb2},  they induce \mlsps s   $\Phi^{(L)}$  on $\Omega \times L$   for $L=\R^N\times L_p([-1,0],\R^N)$,  $1/p+1/q=1$,   and $\Phi^{(C)}$ on  $\Omega \times C$ for $C=C([-1,0],\R^N)$.
Following ideas from~\cite{MiShPart3}, 
 appropriate cooperativity and irreducibility conditions on the systems are assumed to show that both  $\Phi^{(L)}$ and  $\Phi^{(C)}$ admit a family of generalized Floquet subspaces with principal Lyapunov exponent $\widetilde{\lambda}_1^{(L)}$  and $\widetilde{\lambda}_1^{(C)}$ respectively.  In particular, the focusing condition needed for the exponential separation of type II, allows the irreducibility condition to be expressed in terms of the matrix  $A+B$, providing a simplified and weaker version of previous conditions of this type in terms of the matrix $B$.
Since $L$ is separable,  an application of Theorem~\ref{thm:exponential-separation} provides a generalized exponential separation of type II for   $\Phi^{(L)}$  with $\widetilde{\lambda}_1^{(L)}>-\infty$. 
 This result cannot be obtained for $\Phi^C$  in this way  because $C$ is not separable. However, the natural injection $J\colon C\rightarrow L$, introduced and discussed in~\cite{MiNoOb2}, which in particular proves the norm-equivalence of the trajectories defined in both Banach spaces,  allows  us to deduce the result for $\Phi^C$ in Theorem~\ref{expsepC}. For similar reasoning in the context of some fluid mechanics equations see Blumenthal and Punshon-Smith~\cite{Blu-PS}.
 \par\smallskip
Finally, assuming that $\Omega$ is a Lebesgue space and by means of the results from the Oseledets theory  obtained  in Section~\ref{sect:Oseledets},  Section~\ref{sect:lastsection}  shows  the existence of a  exponential separation of time II  in some illustrative examples where the Banach space is not separable.
 More precisely, the above family of linear random delay systems,  under two different $L_1$-integrability conditions  on $B$, induce   \mlsps s   $\Phi^{(\widehat{L})}$  on $\Omega \times \widehat{L}$ for the separable Banach space $\widehat{L}=\R^N\times L_1([-1,0],\R^N)$  and  $\Phi^{(C)}$  on $\Omega \times C$ for $C=C([-1,0],\R^N)$. Then, with the same methods and assumptions of cooperativity and irreducibility of Section~\ref{sect:semiflows},  both  $\Phi^{(\widehat L)}$ and  $\Phi^{(C)}$ admit a family of generalized Floquet subspaces. Lastly, conclusions of Section~\ref{sect:Oseledets} provide a generalized exponential separation of type~II for them with principal Lyapunov exponent $\widetilde{\lambda}_1^{(\widehat{L})}>-\infty$ and $\widetilde{\lambda}_1^{(C)}>-\infty$. The third case corresponds to the separable Banach space of absolutely continuous functions $AC=AC([-1,0],\RR^N)$  with a Sobolev type norm which is not monotone and the cone is not normal, to show the applicability under a weaker condition on the positive cone.  First we consider the same assumptions of the previous case for $\Phi^{(C)}$  to obtain a \mlsps\ $\Phi^{(AC)}$. Then, from the results for $C$, Section~\ref{sect:Oseledets} and similar arguments to those of Theorem~\ref{expsepC}, a family of generalized  Floquet subspaces  and a  exponential separation of type II are obtained for $AC$.
\section{A direct theory providing  generalized exponential separation}\label{sec-prel}
A {\em probability space\/} is a triple $\OFP$, where $\Omega$ is a set, $\mathfrak{F}$ is a $\sigma$\nobreakdash-\hspace{0pt}algebra of subsets of $\Omega$, and $\PP$ is a probability measure defined for all $F \in \mathfrak{F}$.  We always assume that the measure $\PP$ is complete. For a metric space $S$ by $\mathfrak{B}(S)$ we denote the $\sigma$\nobreakdash-\hspace{0pt}algebra of Borel subsets of $S$.
\par\smallskip
A \emph{measurable dynamical system} on the probability space $\OFP$ is a $(\mathfrak{B}(\R) \otimes \mathfrak{F},\mathfrak{F})$\nobreakdash-\hspace{0pt}measurable mapping $\theta\colon\R\times \Omega\to \Omega$ such that
\begin{itemize}
\item
$\theta(0,\w)=\w$ for any $\w\in\Omega$,
\item
$\theta(t_1+t_2,w)=\theta(t_2,\theta(t_1,\w))$ for any $t_1$, $t_2\in\R$ and any $\w \in\Omega$.
\end{itemize}
We write $\theta(t,\w)$ as $\theta_t\w$. Also, we usually denote measurable dynamical systems by $(\OFP,(\theta_{t})_{t \in \R})$ or simply by $(\theta_{t})_{t \in \R}$.\par
A \emph{metric dynamical system} is a measurable dynamical system $(\OFP,(\theta_{t})_{t \in \R})$
such that for each $t\in\R$  the mapping $\theta_t\colon \Omega\to\Omega$ is $\PP$-preserving (i.e., $\PP(\theta_t^{-1}(F))=\PP(F)$ for any $F\in\mathfrak{F}$ and $t\in\R$).
A subset $\Omega'\subset\Omega$ is \emph{invariant} if $\theta_t(\Omega')=\Omega'$ for all $t\in\R$, and the metric dynamical system is said to be \emph{ergodic} if for any invariant subset  $F \in \mathfrak{F}$, either $\PP(F) = 1$ or $\PP(F) = 0$. Throughout the paper we will assume that $\PP$ is ergodic.
\subsection{Measurable linear skew-product semidynamical systems}\label{subsec-measurable}
We consider a separable Banach space $X$ with dual space $X^{*}$.
\par\smallskip
We write $\R^{+}$ for $[0, \infty)$.  By a {\em measurable linear skew-product semidynamical system or semiflow},
$\Phi = \allowbreak ((U_\w(t))_{\w \in \Omega, t \in \R^{+}}, \allowbreak (\theta_t)_{t\in\R})$ on  $X$ covering a metric dynamical system $(\theta_{t})_{t \in \R}$ we understand a $(\mathfrak{B}(\R^{+}) \otimes \mathfrak{F} \otimes \mathfrak{B}(X), \mathfrak{B}(X))$\nobreakdash-\hspace{0pt}measurable
mapping
\begin{equation*}
[\, \R^{+} \times \Omega \times X \ni (t,\w,u) \mapsto U_{\w}(t)\,u \in X \,]
\end{equation*}
satisfying
\begin{align}
&U_{\w}(0) = \mathrm{Id}_{X} \quad & \textrm{for each }\,\w  \in \Omega, \nonumber 
\\
&U_{\theta_{s}\w}(t) \circ U_{\w}(s)= U_{\w}(t+s) \qquad &\textrm{for each } \,\w \in \Omega \textrm{ and }  t,\,s \in \R^{+},
\label{eq-cocycle}
\\
&[\, X \ni u \mapsto U_{\w}(t)u \in X \,] \in \mathcal{L}(X) & \textrm{for each }\,\w \in \Omega \textrm{ and } t \in \R^{+}.\nonumber
\end{align}
Sometimes we write simply $\Phi = ((U_\w(t)), \allowbreak (\theta_t))$. Eq.~\eqref{eq-cocycle} is called the {\em cocycle property\/}.
\par\smallskip
For $\w \in \Omega$, by an {\em entire orbit\/} of $U_{\w}$ we understand a mapping $v_{\w} \colon \R \to X$ such that $v_\w(s + t) = U_{\theta_{s}\w}(t)\, v_\w(s)$ for each $s \in \R$ and $t\geq 0$.
\par\smallskip
Next we introduce the  {\em dual\/} of $\Phi$ in the case in which $X^{*}$ is separable. $\langle \cdot, \cdot \rangle$ will stand for the duality pairing: $\langle u, u^{*} \rangle$ is the action of a functional $u^{*} \in X^{*}$ on a vector $u \in X$.  For $\w \in \Omega$, $t \in \R^{+}$ and $u^{*} \in X^*$ we define $U^{*}_{\w}(t)\,u^{*}$~by
\begin{equation*}
\langle u, U^{*}_{\w}(t)\,u^{*} \rangle = \langle U_{\theta_{-t}\w}(t)\,u , u^{*} \rangle \qquad \text{for each } u \in X
\end{equation*}
(in other words, $U^{*}_{\w}(t)$ is the mapping dual to $U_{\theta_{-t}\w}(t)$).
\par\smallskip
As explained in~\cite{MiShPart1}, since $X^*$ is separable, the mapping
\begin{equation*}
[\, \R^{+} \times \Omega \times X^{*} \ni (t,\w,u^{*}) \mapsto U^{*}_{\w}(t)\,u^{*} \in X^{*} \,]
\end{equation*}
is $(\mathfrak{B}(\R^{+}) \otimes \mathfrak{F} \otimes \mathfrak{B}(X^{*}),
\mathfrak{B}(X^{*}))$\nobreakdash-\hspace{0pt}measurable. The \mlsps\ $\Phi^{*} = ((U^{*}_\w(t))_{\w \in \Omega, t \in \R^{+}},(\theta_{-t})_{t\in\R})$ on $X^{*}$ covering $(\theta_{-t})_{t \in \R}$ will be called  the {\em dual\/} of $\Phi$. The cocycle property for the dual takes the form
\begin{equation*}
U^{*}_{\theta_{-t}\w}(s) \circ U^{*}_{\w}(t)= U^{*}_{\w}(t+s) \qquad \textrm{for each } \,\w \in \Omega \textrm{ and }  t,\,s \in \R^{+}\,.
\end{equation*}
\par\smallskip
Let $\Omega_0 \in \mathfrak{F}$.
A family $\{E(\w)\}_{\w \in \Omega_0}$ of $l$\nobreakdash-\hspace{0pt}dimensional vector subspaces of $X$ is {\em measurable\/} if there are $(\mathfrak{F},
\mathfrak{B}(X))$-measurable functions $v_1, \dots, v_l \colon \Omega_0 \to X$ such that $(v_1(\w), \dots, v_l(\w))$ forms a basis of $E(\w)$ for each $\w \in \Omega_0$.
\par\smallskip
Let $\{E(\w)\}_{\w \in \Omega_0}$ be a family of $l$\nobreakdash-\hspace{0pt}dimensional vector subspaces of $X$, and let $\{F(\w)\}_{\w \in \Omega_0}$ be a family of $l$\nobreakdash-\hspace{0pt}codimensional closed vector subspaces of $X$, such that $E(\w) \oplus F(\w) = X$ for all $\w \in \Omega_0$.  We define the {\em family of projections associated with the decomposition\/} $E(\w) \oplus F(\w) = X$ as $\{P(\w)\}_{\w \in \Omega_0}$, where $P(\w)$ is the linear projection of $X$ onto $F(\w)$ along $E(\w)$, for each $\w \in \Omega_0$.
\par\smallskip
The family of projections associated with the decomposition $E(\w) \oplus F(\w) = X$ is called {\em strongly measurable\/} if for each $u \in X$ the mapping $[\, \Omega_0 \ni \w \mapsto P(\w)u \in X \,]$ is $(\mathfrak{F}, \mathfrak{B}(X))$\nobreakdash-\hspace{0pt}measurable.
\par\smallskip
We say that the decomposition $E(\w) \oplus F(\w) = X$, with
$\{E(\w)\}_{\w \in \Omega_0}$ finite\nobreakdash-\hspace{0pt}dimensional, is {\em invariant\/} if $\Omega_0$ is invariant, $U_{\w}(t)E(\w) = E(\theta_{t}\w)$
and $U_{\w}(t)F(\w) \subset F(\theta_{t}\w)$, for each $t \in \R^+$.
\par\smallskip
A strongly measurable family of projections associated with the invariant decomposition $E(\w) \oplus F(\w) = X$ is referred to as {\em tempered\/} if
\begin{equation*} 
\lim\limits_{t \to \pm\infty} \frac{\ln{\n{P(\theta_{t}\w)}}}{t} = 0 \qquad \PP\text{-a.e. on }\Omega_0.
\end{equation*}
\begin{nota}
\label{rem:tempered-definition}
    In the present paper, instead of the family $\{P(\w)\}_{\w \in \Omega_0}$ of projections  onto $F(\w)$ along $E(\w)$ we usually employ the family $\{\widetilde{P}(\w)\}_{\w \in \Omega_0}$ of projections  onto $E(\w)$ along $F(\w)$ which are related by $\widetilde{P}(\w) = \mathrm{Id}_X - P(\w)$.  It is straightforward that  for the definition of strong measurability we can check the $(\mathfrak{F}, \mathfrak{B}(X))$\nobreakdash-\hspace{0pt}measurability of the mapping $[\, \Omega_0 \ni \w \mapsto \widetilde{P}(\w) \, u \in X \,]$ for each $u \in X$.  Similarly, due the inequalities $1 \le \norm{\widetilde{P}(\w) } \le 1 + \norm{P(\w) }$   and $1 \le \norm{P(\w) } \le 1 + \norm{\widetilde{P}(\w) }$,
    a strongly measurable family of projections associated is tempered if~and only~if
    \begin{equation*}
        \lim\limits_{t \to \pm \infty} \frac{\ln\norm{\widetilde{P}(\theta_{t}\w)}}{t} = 0 \qquad \mathbb{P}\text{-a.e.\ on }\Omega_0.
    \end{equation*}
\end{nota}
\subsection{Ordered Banach spaces}\label{subsec-OBS} Let $X$ be a Banach space with norm $\n{\cdot}$. We say that $X$ is an {\em ordered Banach space\/} if there is a closed convex cone, that is, a nonempty closed subset $X^+\subset X$ satisfying
\begin{itemize}
\item[(O1)] $X^++X^+\subset X^+$.
\item[(O2)] $\R^+ X^+\subset X^+$.
\item[(O3)] $X^+\cap (-X^+)=\{0\}$.
\end{itemize}
Then  a partial ordering in $X$ is defined by
\begin{align*}
u  \leq v \quad& \Longleftrightarrow \quad v -u\in X_+\,;\\
 u< v    \quad  & \Longleftrightarrow \quad v-u\in X_+
 \text{ and }\; u\neq v\,.
\end{align*}
The cone $X^+$ is said to be {\em reproducing\/} if $X^+ - X^+ = X$. The cone $X^+$ is said to be {\em normal} if the norm of the Banach space $X$ is {\em semimonotone\/}, i.e., there is a positive constant $k>0$ such that $0\leq u\leq v$ implies $\n{u}\leq k\,\n{v}$. In such a case, the Banach space can be renormed so that for any $u$, $v\in X$, $0\leq u\leq v$ implies $\n{u}\leq\n{v}$ (see Schaefer~\cite[V.3.1, p. 216]{Schaef}). Such a  norm  is called {\em monotone}.
\par
\smallskip
For an ordered Banach space $X$ denote by $(X^{*})^{+}$ the set of all $u^{*} \in X^{*}$ such that $\langle u, u^{*} \rangle \ge 0$ for all $u \in X^{+}$.  The set $(X^{*})^{+}$ has the properties of a cone, except that $(X^{*})^{+} \cap (-(X^{*})^{+}) = \{0\}$ need not be satisfied (such sets are called {\em wedges\/}).
\par
\smallskip
If $(X^{*})^{+}$ is a cone we call it the {\em dual cone}.  This happens, for instance, when $X^{+}$ is total (that is, $X^{+}- X^{+}$ is dense in $X$, which in particular holds when $X^+$ is reproducing and this will be one of our hypothesis).
Nonzero elements of $X^{+}$ (resp.~of $(X^{*})^{+}$) are called {\em positive\/}.
\subsection{Assumptions}\label{sec-assumptions} Throughout the paper, we will assume that $X$ is an ordered separable  Banach space with $\dim X\geq 2$, dual space $X^{*}$ and  positive cone $X^+$ normal and reproducing.
\par
\smallskip
Let $\Phi = ((U_\w(t)), (\theta_t))$ be a \mlsps\ on $X$ covering an ergodic metric dynamical system $(\theta_{t})$ on $\OFP$. We now list assumptions we will make at various points in the sequel.
\par\smallskip
\begin{enumerate}[label=(\textbf{A\arabic*}),series=Assumptions,leftmargin=28pt]
\item\label{A1} (Integrability) The functions
\begin{equation*}
\bigl[ \, \Omega \ni \w \mapsto \sup\limits_{0 \le s \le 1} {\lnplus{\n{U_{\w}(s)}}} \in [0,\infty) \, \bigr] \in L_1\OFP\;\; \text{ and}
\end{equation*}
\begin{equation*}
\bigl[ \, \Omega \ni \w \mapsto \sup\limits_{0 \le s \le 1}
{\lnplus{\n{U_{\theta_{s}\w}(1-s)}}} \in [0,\infty) \, \bigr] \in L_1\OFP.
\end{equation*}
\end{enumerate}
\begin{enumerate}[resume*=Assumptions]
\item\label{A2} (Positivity) For any $\w \in \Omega$, $t \ge 0$ and $u_1, u_2 \in X$ with $u_1 \le u_2$
\begin{equation*}
U_{\w}(t)\,u_1 \le U_{\w}(t)\,u_2\,.
\end{equation*}
\end{enumerate}
\par
\smallskip
For a \mlsps\ $\Phi$ satisfying  \ref{A2} we say that an entire orbit $v_{\w}$ of $U_{\w}$ is {\em positive\/} if $v_{\w}(t) \in X^{+} \setminus \{0\}$ for all $t \in \R$.
\par
\smallskip
Next we introduce focusing condition \ref{A3} in the following way.
\par\smallskip
\begin{enumerate}[resume*=Assumptions]
\item\label{A3} (Focusing) \ref{A2} is satisfied  and there are $\be\in X^+$ with $\|\be\|=1$ and $U_\w(1)\,\be\neq 0 $ for all $\w\in\Omega$, and an $(\mathfrak{F},\mathfrak{B}(\R))$\nbd-measurable function $\varkappa\colon\Omega\to [1,\infty)$ with
$\lnplus\ln\varkappa\in L_1(\Omega,\mathfrak{F},\PP)$ such that for any  $\w\in\Omega$ and any nonzero $u\in X^+$
\begin{itemize}
\item $U_\w(1)\,u=0$, or
\item $U_\w(1)\,u\neq 0$ and there is $\beta(\w,u)>0$ with the property that
\begin{equation*}
\beta(\w,u)\,\be\leq U_\w(1)\,u\leq \varkappa(\w)\,\beta(\w,u)\,\be\,.
\end{equation*}
\end{itemize}
\end{enumerate}
It should be remarked that our condition \ref{A3} is weaker than condition \ref{A3} in~\cite{MiShPart1}.
\begin{nota}\label{rm:dichotomy}
Under~\ref{A3}
for $u \in X^{+}$ the following dichotomy holds:
\begin{itemize}
\item
$U_\w(t)\, u = 0$ for all $t \ge 1$, or
\item
$U_\w(t)\, u > 0$ for all $t \ge 1$.
\end{itemize}
\end{nota}
When the dual space $X^*$ is separable, we consider the dual of $\Phi$, as explained in Subsection~\ref{subsec-measurable}, $\Phi^{*} = ((U^{*}_\w(t))_{\w \in \Omega, t \in \R^{+}},(\theta_{-t})_{t\in\R})$ on $X^{*}$ covering $(\theta_{-t})_{t \in \R}$, and the following assumptions for it.
\begin{enumerate}[label=(\textbf{A\arabic*})*,series=Assumptions*,leftmargin=32pt]
\item\label{A1*} \!\!(Integrability) The functions
\begin{equation*}
\bigl[ \, \Omega \ni \w \mapsto \sup\limits_{0 \le s \le 1} {\lnplus{\n{U_{\w}^*(s)}}} \in [0,\infty) \, \bigr] \in L_1\OFP\;\; \text{ and}
\end{equation*}
\begin{equation*}
\bigl[ \, \Omega \ni \w \mapsto \sup\limits_{0 \le s \le 1}
{\lnplus{\n{U_{\theta_{s}\w}^*(1-s)}}} \in [0,\infty) \, \bigr] \in L_1\OFP.
\end{equation*}
\end{enumerate}
\begin{enumerate}[resume*=Assumptions*]
\item\label{A2*} \!\!(Positivity)  For any $\w \in \Omega$, $t \ge 0$ and $u^{*}_1, u^{*}_2 \in X^*$ with $u^{*}_1 \le u^{*}_2$
\begin{equation*}
U^{*}_{\w}(t)\,u^{*}_1 \le U^{*}_{\w}(t)\,u^{*}_2\,.
\end{equation*}
\end{enumerate}
Notice that in this case, as explained in~\cite{MiShPart1}, \ref{A1*} follows from \ref{A1} and \ref{A2*} follows from \ref{A2}.
\begin{enumerate}[resume*=Assumptions*]
\item\label{A3*}  (Focusing for $X^*$) \ref{A2*} is satisfied and there are ${\bf e^*}\in (X^*)^+$ with $\langle \be,\be^*\rangle=1$ and $\| \be^*\|=1$ and an $(\mathfrak{F},\mathfrak{B}(\R))$-measurable function $\varkappa^*\colon \Omega\to[1,\infty)$ with
$\lnplus\ln\varkappa^*\in L_1(\Omega,\mathfrak{F},\PP)$ such that for any  $\w\in\Omega$ there holds $U_\w^*(1)\, \be^* \ne 0$, and for any $\w\in\Omega$ and any nonzero $u^*\in (X^*)^+$ there is $\beta^*(\w,u^*)>0$ with the property that
\begin{equation*}
\beta^*(\w,u)\,U_\w^*(1)\,{\bf e^*}\leq U_\w^*(1)\,u^*\leq \varkappa^*(\w)\,\beta^*(\w,u^*)\,U^*_\w(1)\,\be^*\,.
\end{equation*}
\end{enumerate}
Condition \ref{A3*} follows from \ref{A3}, as proved in~\cite[Proposition 2.4]{MiNoOb1},
\begin{nota}\label{timeT}
We can replace time 1 with some nonzero $T\in\R^+$ in~\textup{\ref{A1}}, \ref{A3}, and  \ref{A1*}, \ref{A3*}.
\end{nota}
\subsection{Generalized exponential separation}\label{sec-GES}
 In this subsection we recall the definition and existence of a family of generalized principal Floquet subspaces, and of a generalized exponential separation of type II, introduced and proved in~\cite{MiNoOb1}. They are important for the cases in which the previous concepts given in~\cite{MiShPart1} do not apply, as \mlsps s  induced by delay differential equations.
\begin{defi}
\label{generalized-floquet-space}
Let $\Phi = ((U_\w(t)), (\theta_t))$ be a measurable linear skew-product semidynamical system on a Banach space $X$ ordered by a normal reproducing cone $X^{+}$.  A family of one\nobreakdash-\hspace{0pt}dimensional subspaces $\{E_1(\w)\}_{\w \in \widetilde{\Omega}}$ of $X$ is called a family of {\em generalized principal Floquet subspaces} of $\Phi$ if $\widetilde{\Omega} \subset \Omega$ is invariant, $\PP(\widetilde{\Omega}) = 1$, and
\begin{itemize}
\item[{\rm (i)}]
$E_1(\w) = \spanned{\{w(\w)\}}$ with $w \colon \widetilde{\Omega} \to X^+ \setminus \{0\}$ being $(\mathfrak{F}, \mathfrak{B}(X))$\nobreakdash-\hspace{0pt}measurable,
\item[{\rm (ii)}]
$U_{\w}(t) \,E_1(\w) = E_1(\theta_{t}\w)$, for any $\w \in \widetilde{\Omega}$ and any $t > 0$,
\item[\textup{(iii)}]
there exists $\widetilde{\lambda} \in [-\infty, \infty)$ such that
\begin{equation*}
\widetilde{\lambda} = \lim_{t\to\infty} \frac{\ln{\n{U_\w(t)\,w(\w)}}}{t} \quad \text{ for each } \w \in \widetilde{\Omega},
\end{equation*}
\item[{\rm (iv)}] for each $\w\in\widetilde\Omega$ and each $ u\in X^+$ with $U_\w(1)\,u\neq 0$
\begin{equation}\label{principal-exponent}
\limsup_{t\to\infty} \frac{\ln{\n{U_\w(t)\,u}}}{t} \le \widetilde{\lambda} \,.
\end{equation}
\end{itemize}
$\widetilde{\lambda}$ is called the {\em generalized principal Lyapunov exponent} of $\Phi$ associated to the generalized principal Floquet subspaces $\{E_1(\w)\}_{\w\in\widetilde{\Omega}}$.
\end{defi}
Under assumptions {\rm \ref{A1}} and {\rm\ref{A3}}, (i)--(v) of the next theorem, proved in~\cite[Theorem 3.10]{MiNoOb1}, shows the existence of an invariant set $\widetilde\Omega_1$ of full measure $\PP(\widetilde\Omega_1)=1$, a family of generalized Floquet subspaces  $\{E_1(\w)\}_{\w \in \widetilde\Omega_1}$, with $E_1(\w) = \spanned\{w(\w)\}$, and generalized principal Lyapunov exponent  $\widetilde{\lambda}_1$.
Notice that  in~\cite{MiNoOb1} a standing assumption is that not only $X$ but also its dual $X^{*}$ is a separable Banach space.  However, the next result consider \mlsps{}s taking values in $X$ only, so no assumptions on $X^{*}$ are needed.
\begin{teor}\label{Theorem3.10}
Under assumptions   {\rm \ref{A1}} and {\rm\ref{A3}}, there is an invariant set  $\bar\Omega_1$ and an $(\mathfrak{F},\mathfrak{B}(X))$-measurable function
$w\colon\bar\Omega_1\to X^+,\;\w\mapsto w(\w)$,     $\n{w(\w)}=1$ for each $\w\in \bar\Omega_1$,   such that
\begin{itemize}
\item[{\rm (i)}] for each $\w\in\bar\Omega_1$ and $t\geq 0$,
\begin{equation*}
w(\theta_t\w)=\frac{U_\w(t)\,w(\w)}{\n{U_\w(t)\,w(\w)}}\,;
\end{equation*}
\item[{\rm (ii)}]
for each $\w\in\bar\Omega_1$, the map $w_\w\colon \R\to X^+$ defined by
\begin{equation*}
w_\w(t)
= \begin{cases}
\displaystyle\frac{w(\theta_t\w)}{\n{U_{\theta_t\w}(-t) \, w(\theta_t\w)}} & \quad \text{for } t\leq 0\,,
\\[.4cm]
\;U_\w(t)\,w(\w)& \quad \text{for } t\geq 0\,; \end{cases}
\end{equation*}
is a positive entire orbit of $U_\w$, unique up to multiplication by a positive scalar;
\item[{\rm (iii)}]
there are an invariant set $\widetilde\Omega_1 \subset \bar\Omega_1$ with $\PP(\widetilde\Omega_1) = 1$ and  a $\widetilde\lambda_1 \in [-\infty,\infty)$ such that
\begin{equation*}
\widetilde\lambda_1 = \lim_{t\to\pm\infty} \frac{1}{t}\ln\n{w_\w(t)} = \int_\Omega \ln\n{w_{\w'}(1)}\,d\PP(\w')
\end{equation*}
for each $\w \in \widetilde\Omega_1$;
\item[{\rm (iv)}]
for each $\w \in \widetilde\Omega_1$ and $u \in X^+$ with $U_\w(1) \,u \neq 0$
\begin{equation*}
\lim_{t\to\infty} \frac{1}{t}\ln\n{U_\w(t)\,u} = \widetilde\lambda_1\,;
\end{equation*}
\item[{\rm (v)}]
for each $\w \in \widetilde\Omega_1$ and $u \in X$
\begin{equation*}
\limsup_{t\to\infty}\frac{1}{t}\ln\n{U_\w(t)\,u}\leq \widetilde\lambda_1\,;
\end{equation*}
that is,  $\{E_1(\w)\}_{\w \in \widetilde\Omega_1}$, with $E_1(\w) = \spanned\{w(\w)\}$, is a family of generalized principal Floquet subspaces, and $\widetilde\lambda_1$ is the generalized principal Lyapunov exponent;
\item[\textup{(vi)}]
for each $\w \in \widetilde\Omega_1$,
\begin{equation*}
\lim_{t\to\infty}
\frac{1}{t} \ln{\sup \{\, \norm{U_{\w}(t) \, u} :  u \in X^{+} ,\,\n{u}\leq 1}\} = \widetilde{\lambda} _1\,;
\end{equation*}
\item[\textup{(vii)}]
there exists $\sigma > 0$ such that, for any $\w \in \widetilde{\Omega}_1$,
\begin{align*}
\quad\;\limsup\limits_{t \to \infty} \frac{1}{t} \ln\Bigl(\sup{\Bigl\{\, \Bigl\lVert \frac{U_{\theta_{-t}\w}(t) \, u}{\norm{U_{\theta_{-t}\w}(t) \, u}} - w(\w) \Bigr\rVert :  u \in X^{+}, \, U_{\theta_{-t}\w}(1) \, u \ne 0} \, \Bigr\} \Bigr) & \le -\sigma,
\\
\limsup\limits_{t \to \infty} \frac{1}{t} \ln\Bigl(\sup{\Bigl\{\, \Bigl\lVert \frac{U_{\w}(t) \, u}{\norm{U_{\w}(t) \, u}} - w(\theta_{t}\w) \Bigr\rVert : u \in X^{+},\,  U_{\w}(1) \, u \ne 0} \, \Bigr\} \Bigr) & \le -\sigma.
\end{align*}
\end{itemize}
\end{teor}
\begin{proof}
    Parts (i)--(v) are \cite[Theorem 3.10]{MiNoOb1}.  To prove part (vi), notice that in the proof of~\cite[Theorem.~3.10]{MiNoOb1} on p.~6175 the inequality
    \begin{equation*}
        \norm{U_{\w}(t) \, u} \le \frac{\varkappa(\w) \beta(\w, u)}{\beta(\w, w(\w))} \,\norm{U_{\w}(t) \, w(\w)}, \quad t \ge 1,
    \end{equation*}
    holds for all $\w \in \widetilde{\Omega}_1$, where $\varkappa(\cdot)$ and $\beta(\cdot, \cdot)$ are as in~\ref{A3}.  For a fixed $\w \in \widetilde{\Omega}_1$, since $U_{\w}(1)$ is linear bounded, $\sup\{\, \beta(\w, u) : u \in X^{+},\,U_{\w}(1) \, u \ne 0,\, \n{u}\leq 1 \,\} < \infty$.  Consequently,
    \begin{multline*}
        \limsup_{t\to\infty} \frac{1}{t} \ln{\sup \{\, \norm{U_{\w}(t)\,u} :  u \in X^{+} \,\}}
        \\
        = \limsup_{t\to\infty} \frac{1}{t} \ln{\sup \{\, \norm{U_{\w}(t)\,u} :  u \in X^{+},\, \n{u}\leq 1}\,\} \le \widetilde{\lambda}_1 \,.
    \end{multline*}
    The other inequality is straightforward by~Theorem~\ref{Theorem3.10}(iv), and finally,  part (vii) follows along the lines of the proof of~\cite[Proposition~5.5]{MiShPart1} (for a fairly similar reasoning, see p.~6173 below~\cite[Proposition~3.9]{MiNoOb1}).
\end{proof}
\begin{defi}
\label{def:exponential-separation}
Let $\Phi = ((U_\w(t)), (\theta_t))$ be a measurable linear skew-product semidynamical system on a Banach space $X$ ordered by a normal reproducing cone $X^{+}$.  $\Phi$ is said to admit a {\em generalized exponential separation of type} II if there are a family of generalized principal Floquet subspaces $\{E_1(\w)\}_{\w \in \widetilde{\Omega}}$,
and a family of one\nbd-codimensional closed vector subspaces $\{F_1(\w)\}_{\w \in \widetilde{\Omega}}$ of $X$,  satisfying
\begin{itemize}
\item[{\rm (i)}]
$F_1(\w)\cap X^+=\{u\in X^+ : U_\w(1)\,u=0\}$,
\item[{\rm (ii)}]
$X=E_1(\w)\oplus F_1(\w)$ for any $\w\in\widetilde{\Omega}$, where the decomposition is invariant, and the family of projections associated with the decomposition is strongly measurable and tempered,
\item[{\rm (iii)}]
there exists $\widetilde{\sigma}\in (0,\infty]$ such that
\[ \lim_{t\to\infty}\frac{1}{t}\ln \frac{\n{U_\w(t)|_{F_1(\w)}}}{\n{U_\w(t)\,w(\w)}}=-\widetilde{\sigma}
\]
for each $\w\in\widetilde{\Omega}$.
\end{itemize}
It is said that $\{E_1(\cdot),F_1(\cdot),\widetilde{\sigma}\}$ {\em generates a generalized exponential separation of type}~II.
\end{defi}
As stated in the introduction, sometimes a less general concept is used, namely that of generalized exponential separation of type~I.  Then in Definition~\ref{def:exponential-separation}(i) we  would have just $F_1(\w) \cap X^+ = \{0\}$.  Note that the only difference is that now  $F_1(\w)$  may contain those positive vectors $u>0$ for which $U_\w(1)\,u=0$.
For a theory of generalized exponential separation of type~I (called there generalized exponential separation) see the series of papers~\cite{MiShPart1, MiShPart2, MiShPart3}.
\par\smallskip
The next theorem, proved in~\cite[Theorem 4.6]{MiNoOb1} and included here for completeness, shows the existence of a generalized exponential separation of type II. We maintain the notation of~Theorem~\ref{Theorem3.10}.
\begin{teor}
\label{thm:exponential-separation}
Assume that $X^*$ is separable. Under assumptions {\rm \ref{A1}} and {\rm\ref{A3}}, let $\widetilde\lambda_1$ be the generalized principal Lyapunov exponent and assume that $\widetilde\lambda_1 > -\infty$.  Then there is an invariant set $\widetilde\Omega_0\subset \widetilde \Omega_1$ of full measure $\PP(\widetilde\Omega_0)=1$ such that
\begin{enumerate}
\item[{\rm (i)}]
The family $\{ P(\w)\}_{\w\in\widetilde\Omega_0}$ of projections associated with invariant decomposition $ E_1(\w)\oplus  F_1(\w)=X$ is strongly measurable and tempered.
\item[{\rm (ii)}]
$F_1(\w)\cap X^+=\{u\in X^+ : U_\w(1)\,u=0\}$ for any $\w\in\widetilde\Omega_0$.
\item[{\rm (iii)}]
For any $\w\in\widetilde\Omega_0$ and $u\in X\setminus F_1(\w)$ with $U_\w(1)\,u\neq 0$  there holds
 \[ \lim_{t\to\infty}\frac{1}{t}\ln\n{U_\w(t)}=\lim_{t\to\infty}\frac{1}{t}\ln \n{U_\w(t)\,u}=\widetilde{\lambda}_1\,.\]
 \item[{\rm (iv)}] There exists $\widetilde\sigma\in(0,\infty]$ and $\widetilde\lambda_2=\widetilde\lambda_1-\widetilde\sigma$ such that
\[ \lim_{t\to\infty}\frac{1}{t}\ln \frac{\n{U_\w(t)|_{F_1(\w)}}}{\n{U_\w(t)\,w(\w)}}=-\widetilde{\sigma}
\]
and
\[\lim_{t\to\infty}\frac{1}{t}\ln\n{U_\w(t)|_{F_1(\w)}}=\widetilde\lambda_2\]
for each $\w\in\widetilde\Omega_0$.
\end{enumerate}
That is, $\Phi$ admits a generalized exponential separation of type {\rm II}.
\end{teor}
\section{Generalized exponential separation as a consequence of Oseledets decomposition}\label{sect:Oseledets}
Let $\Phi = ((U_\w(t)), (\theta_t))$ be a measurable linear skew-product semidynamical system on a separable Banach space $X$ with $\dim{X} \ge 2$, covering a metric dynamical system $(\theta_{t})$. We always assume that \ref{A1} is satisfied, that is, the functions $\bigl[ \, \Omega \ni \w \mapsto \sup\limits_{0 \le s \le 1} {\lnplus{\n{U_{\w}(s)}}} \in [0,\infty) \, \bigr]$ and $\bigl[ \, \Omega \ni \w \mapsto \sup\limits_{0 \le s \le 1} {\lnplus{\n{U_{\theta_{s}\w}(1-s)}}} \in [0,\infty) \, \bigr]$ belong to $L_1\OFP$.  It follows then from Kingman's subadditive ergodic theorem that there exists $\lambda_{\mathrm{top}} \in [-\infty, \infty)$ such that
\begin{equation*}
\lim\limits_{t \to \infty} \frac{\ln{\norm{U_{\w}(t)}}}{t} = \lambda_{\mathrm{top}}
\end{equation*}
for $\PP$-a.e. $\w \in \Omega$.  $\lambda_{\mathrm{top}}$ is called the  {\em top Lyapunov exponent}.
\par\smallskip
Our starting point is the semi-invertible operator Oseledets\nobreakdash-\hspace{0pt}type theorem \cite[Theorem~3.4]{MiNoOb2}, based on the results in~\cite{GTQu}.  We state here its parts used in the sequel.
\begin{teor}
\label{thm-Oseledets}
Let $\OFP$ be a Lebesgue space.  Assume that $U_{\w}(t_0)$ is compact for all $\w \in \Omega$ and some $t_0\geq 1$.  Let $\lambda_{\mathrm{top}} > -\infty$.  Then
there exist:
\begin{itemize}
\item
an invariant $\Omega' \subset \Omega$, with $\PP(\Omega') = 1$,
\item
an invariant measurable decomposition   $X = E(\w) \oplus F(\w)$, $\w \in \Omega'$, such that the family of projections associated with it is tempered,
\item
$\lambda_2 \in [ - \infty, \lambda_{\mathrm{top}})$
\end{itemize}
with the properties that
\begin{itemize}
\item[\textup{(i)}]
$E(\w)$ has constant finite dimension, $l$, on $\Omega'$,
\item[\textup{(ii)}]
for each $\w \in \Omega'$ and each $t > 0$, $U_{\w}(t)|_{E(\w)}$ is a linear automorphism onto $E(\theta_{t}\w)$,
\item[\textup{(iii)}] for any $\w \in \Omega'$
\begin{equation*}
\lim_{t \to \infty} \frac{\ln{\norm{U_{\w}(t)|_{E(\w)}}}}{t} = \lim_{t \to \infty} \frac{\ln{\norm{(U_{\w}(t)|_{E(\w)})^{-1}}^{-1}}}{t} = \lambda_{\mathrm{top}}\,,
\end{equation*}
\item[\textup{(iv)}] for any $\w \in \Omega'$
\begin{equation*}
\lim_{t \to \infty} \frac{\ln{\norm{U_{\w}(t)|_{F(\w)}}}}{t} = \lambda_2\,.
\end{equation*}
\end{itemize}
The measurable invariant decomposition satisfying the above is unique.
\end{teor}
When $\dim{X} < \infty$ it is possible that $E(\w) = X$ for all $\w \in \Omega'$.  In that case, (iv) is vacuous.
\par
Below, we apply some ideas from Blumenthal's paper~\cite{Blum}, elaborated later in Varzaneh and Riedel~\cite{Varz-Ried}. As in~\cite[Subsection~4.0.4]{Lian-Lu}, define the $l$\nobreakdash-\hspace{0pt}dimensional volume function $\vol \colon X^l \to \RR$ by the formula
\begin{equation}
\label{def:vol}
    \vol(v_1, \ldots, v_l) := \norm{v_l} \, \prod\limits_{i = 1}^{l - 1} \dist(v_i, \spanned{\{v_{i + 1}. \ldots, v_l\}}), \quad v_1, \ldots, v_l \in X.
\end{equation}
\begin{prop}
\label{prop:volume}
    Under the assumptions and notations of \textup{Theorem~\ref{thm-Oseledets}}, for any $\w \in \Omega'$ and any basis $\{v_1, \ldots, v_l\}$ of $E(\w)$ there holds
    \begin{equation*}
        \lim_{n \to \infty} \frac{\ln{\vol(U_{\w}(n) \, v_1, \ldots, U_{\w}(n) \, v_l)}}{n} = l \, \lambda_{\mathrm{top}}.
    \end{equation*}
\end{prop}
\begin{proof}   See \cite[Thm.~1.21(v)]{Varz-Ried}.
\end{proof}
The following result is known as one of the Kre\u{\i}n--Shmulyan theorems (see~\cite[Theorem~2.2]{AbAlBu}).
\begin{lema}
\label{lemma:Krein-Shmulyan}
Let $X^{+}$ be a  reproducing cone in a Banach space $X$ with
norm  $\norm{\cdot}$.
Then
there exists $K \ge 1$ with the property that for each $u \in X$ there are $u^{+}, u^{-} \in
X^{+}$ such that $u = u^{+} - u^{-}$, $\norm{u^{+}} \le K \norm{u}$ and $\norm{u^{-}} \le K \norm{u}$.
\end{lema}
From now on, all the results in the present section are obtained under the assumptions of Theorems~\ref{Theorem3.10} and~\ref{thm-Oseledets}.
\begin{teor}
\label{thm-nontrivial}
 There holds
  $E(\w) \cap X^{+} \varsupsetneq \{0\}$ $\PP$-a.e.\ on $\Omega$.
\end{teor}
\begin{proof}
It goes along the lines of the proof of Theorem~3.5 in~\cite{MiShPart1}.
\end{proof}
From now on until the end of the present section we assume that, whenever we talk about generalized principal Floquet subspaces, the vectors $w(\w) \in X^{+} \setminus \{0\}$ spanning the one\nobreakdash-\hspace{0pt}dimensional subspace are chosen to be unit vectors.
\begin{lema}  \label{lm:top-equals-principal}
   $\lambda_{\mathrm{top}} = \widetilde{\lambda}_1$.
\end{lema}
\begin{proof}
    The inequality $\widetilde{\lambda} _1\le \lambda_{\mathrm{top}}$ is straightforward.  In order to prove the other one, observe that
\begin{align*}
        \norm{U_{\w}(t)} & = \sup\{\, U_{\w}(t) (u^{+} - u^{-}) : \norm{u^{+} - u^{-}} \le 1 \,\}
        \\
        & \le 2\, K \sup\{\, U_{\w}(t) \, u :  u \in X^{+}, \norm{u} \le 1 \,\},
\end{align*}
    where $K \ge 1$ is a constant in Lemma~\ref{lemma:Krein-Shmulyan}, and apply Theorem~\ref{Theorem3.10}(vi).
\end{proof}
We should mention here that a similar result, \cite[Proposition~2.2]{MiShPart3}, was proved with the help of Baire's theorem.
\begin{lema}
\label{lemma-subspace}
$E_1(\w) \subset E(\w)$ for $\PP$-a.e.\ $\w \in \Omega$.
\end{lema}
\begin{proof}
It is sufficient to show that $w(\w) \in E(\w)$ for $\PP$-a.e.\ $\w \in \Omega$.  By Theorem~\ref{thm-nontrivial}, $\PP$-a.e.\ on $\Omega$, for each $n \in \NN$ we can find $u_n \in E(\theta_{-n}\w) \cap X^{+}$, with $\norm{u_n} = 1$.  The invariance of $E$ gives that $U_{\theta_{-n}\w}(n) \, u_n/\norm{U_{\theta_{-n}\w}(n) \, u_n}$ belongs to $E(\w)$.  By Theorem~\ref{Theorem3.10}(vii), $U_{\theta_{-n}\w}(n) \, u_n/\norm{U_{\theta_{-n}\w}(n) \, u_n}$ converges, as $n \to \infty$, to $w(\w)$.
\end{proof}
\begin{teor}
\label{thm-identity}
The dimension $l$ of $E(\w)$ equals one.
\end{teor}
\begin{proof}
 Let $\widetilde{\Omega}_1$ be the invariant subset of $\Omega$ of Theorem~\ref{Theorem3.10} and without loss of generality, assume that the set $\Omega'$ of Theorem~\ref{thm-Oseledets} is valid for Theorem~\ref{thm-nontrivial}, that is,  $E(\w) \cap X^{+} \varsupsetneq \{0\}$ for all $\w \in \Omega'$. \par\smallskip
Suppose to the contrary that $l \ge 2$.
Fix $\w \in \Omega' \cap \widetilde{\Omega}_1$ such that $E_1(\w) \subset E(\w)$, and let $\{u_1, \ldots, u_{l - 1}, w(\w) \}$ be a basis of unit vectors for $E(\w)$.  We will look at the action of $U_{\w}(n)$ on  $l$\nobreakdash-\hspace{0pt}dimensional volume.
We have, by~\eqref{def:vol},
\begin{multline}\label{eq:volume-1}
\hspace{-.4cm}\vol(U_{\w}(n) \, u_1, \ldots, U_{\w}(n) \, u_{l - 1}, U_{\w}(n) \, w(\w) )
\\ \hspace{.062cm}= \norm{U_{\w}(n) \, w(\w)} \, \prod\limits_{i = 1}^{l - 1} \dist(U_{\w}(n) \, u_i, \spanned{\{U_{\w}(n) \, u_{i + 1}, \ldots, U_{\w}(n) \, u_{l - 1}, U_{\w}(n) \, w(\w)\}})
 \\[-.2cm]
 \le \norm{U_{\w}(n) \, w(\w)} \, \prod\limits_{i = 1}^{l - 1} \dist(U_{\w}(n) \, u_i, E_1(\theta_{n}\w))
\end{multline}
for all $n \in \NN$.
\par\smallskip
By Lemma~\ref{lemma:Krein-Shmulyan}, there is $K \ge 1$ such that for each $u_i$, $1 \le i \le l - 1$, one can find $u_{i}^{+}, u_{i}^{-} \in X^{+}$ with $u_{i} = u_{i}^{+} - u_{i}^{-}$, $\norm{u_{i}^{+}} \le K$, $\norm{u_{i}^{-}} \le K$.
\par\smallskip
It follows from Theorem~\ref{Theorem3.10}(vii) that for any $\rho \in (0, \sigma)$ one can find $N_0$ such that for $n \ge N_0$ and $i \in \{1, \ldots, l - 1 \}$ there holds
\begin{gather*}
\Bigl\lVert U_{\w}(n) \, u_{i}^{+} - \norm{U_{\w}(n) \, u_{i}^{+}} \, w(\theta_{n}\w) \Bigr\rVert \le \exp(- n \rho) \norm{U_{\w}(n) \, u_{i}^{+}},
\\
\Bigl\lVert U_{\w}(n) \, u_{i}^{-} - \norm{U_{\w}(n) \, u_{i}^{-}} \, w(\theta_{n}\w) \Bigr\rVert  \le \exp(- n \rho) \norm{U_{\w}(n) \, u_{i}^{-}},
\end{gather*}
consequently,
\begin{multline}
\label{eq:volume-2}
\Bigl\lVert U_{\w}(n) \, u_{i} - \bigl( \norm{U_{\w}(n) \, u_{i}^{+}} - \norm{U_{\w}(n) \, u_{i}^{-}} \bigr) \, w(\theta_{n}\w) \Bigr\rVert
\\
= \Bigl\lVert \bigl(U_{\w}(n) \, u_{i}^{+} - \norm{U_{\w}(n) \, u_{i}^{+}} \, w(\theta_{n}\w) \bigr) - \bigl( U_{\w}(n) \, u_{i}^{-} - \norm{U_{\w}(n) \, u_{i}^{-}} \, w(\theta_{n}\w) \bigr) \Bigr\rVert
\\
\le \exp(- n \rho) \bigl( \norm{U_{\w}(n) \, u_{i}^{+}} + \norm{U_{\w}(n) \, u_{i}^{-}} \bigr).
\end{multline}
Pick some $\mu$  satisfying
\begin{equation*}
    \lambda_{\mathrm{top}} < \mu < \lambda_{\mathrm{top}} + \frac{l - 1}{l} \,\rho.
\end{equation*}
There is $N_1$ such that for $n \ge N_1$ and $i \in \{1, \ldots, l - 1 \}$ there holds
\begin{equation}
\label{eq:volume-3}
\norm{U_{\w}(n) \, u_{i}^{+}} \le \exp(n \mu) \norm{u_{i}^{+}}, \quad \norm{U_{\w}(n) \, u_{i}^{-}} \le \exp(n \mu) \norm{u_{i}^{-}}.
\end{equation}
Gathering \eqref{eq:volume-2} and \eqref{eq:volume-3} we obtain, in view of Lemma~\ref{lemma:Krein-Shmulyan}, the following:
\begin{equation*}
\dist(U_{\w}(n) \, u_i, E_1(\theta_{n}\w)) \le 2 K \exp(n (\mu - \rho)), \quad 1 \le i \le l - 1,
\end{equation*}
for $n \ge \max\{N_0, N_1\}$, which gives, via~\eqref{eq:volume-1},
\begin{multline*}
\vol(U_{\w}(n) \, u_1, \ldots, U_{\w}(n) \, u_{l - 1}, U_{\w}(n) \, w(\w) )
\\
\le 2^{l-1} K^{l-1} \exp((l-1) n (\mu - \rho)) \, \exp(n \mu),
\end{multline*}
consequently,
\begin{equation*}
\limsup\limits_{n \to \infty} \frac{\ln{\vol(U_{\w}(n) \, u_1, \ldots, U_{\w}(n) \, u_{l - 1}, U_{\w}(n) \, w(\w) )}}{n} \le l \mu - (l - 1) \rho,
\end{equation*}
which contradicts Proposition~\ref{prop:volume} and finishes the proof.
\end{proof}
As a consequence, we can deduce the existence of a generalized exponential separation of type II as shown in the following theorem.
\begin{teor}
\label{thm:Appendix}
Let  $\OFP$ be a Lebesgue space, $\Phi = ((U_\w(t)), (\theta_t))$ be a \mlsps\ on a separable Banach space $X$ with $\dim{X} \ge 2$, positive cone $X^+$ normal and reproducing.   Assume that  \textup{\ref{A1}} and \textup{\ref{A3}} hold,  $U_{\w}(t_0)$ is compact for all $\w \in \Omega$ and some $t_0\geq 1$, and moreover $\lambda_{\mathrm{top}} > -\infty$.  Then $((U_\w(t)), \allowbreak (\theta_t))$ admits a generalized exponential separation of type \textup{II}.
\end{teor}
\begin{proof}
We claim that $F(\w)$   of Theorem~\ref{thm-Oseledets} can serve as $F_1(\w)$ in the definition of generalized exponential separation.
Indeed, parts (ii) and (iii) in Definition~\ref{def:exponential-separation} are direct consequences of Theorems~\ref{thm-Oseledets} and~\ref{Theorem3.10} combined with Theorem~\ref{thm-identity}.
 In order to prove (i), let $u \in X^{+}$ be such that $U_{\w}(1) \, u = 0$.  Then, 
 $\lim_{t \to \infty} \ln{\norm{U_{\w}(t) \, u}}/t = - \infty$,
 and hence, by the characterization given in Theorem~\ref{thm-Oseledets}(iii)--(iv), $u \in F(\w) = F_1(\w)$.  Finally, suppose to the contrary that there is a $u \in F(\w) \cap X^{+}$ such that $U_{\w}(1) \, u \ne 0$.  Then, by Remark~\ref{rm:dichotomy}, $U_{\w}(t) \, u \in X^{+} \setminus \{0\}$ for all $t \ge 1$. Hence, by Theorem~\ref{Theorem3.10}(ii) and Lemma~\ref{lm:top-equals-principal}, $\lim_{t \to \infty} \ln{\norm{U_{\w}(t) \, u}}/t = \lambda_{\mathrm{top}}$, which contradicts Theorem~\ref{thm-Oseledets}(iv),  shows that Definition~\ref{def:exponential-separation}(i) holds and finishes the proof.
\end{proof}
We should mention here that Theorem~\ref{thm:Appendix} is new even in the case of generalized exponential separation of type I.  Indeed, in~\cite[Theorem.~3.8]{MiShPart1} generalized exponential separation of type I was proved (with stipulating neither $\OFP$ to be a Lebesgue space nor $\lambda_{\mathrm{top}} > - \infty$) under a much stronger assumption than (A3), namely that $\ln{\varkappa}\in L_1(\Omega,\mathfrak{F},\PP)$ (such a property was called strong focusing).  As both the Lebesgue property of $\OFP$ and $\lambda_{\mathrm{top}} > - \infty$ are quite natural, the above theorem is, for all conceivable practical purposes, a strengthening of \cite[Theorem~3.8]{MiShPart1}.
\section{Semiflows generated by linear random delay differential equations}\label{sect:semiflows}
This section is devoted to showing the applications of the  theory of Section~\ref{sec-prel} to random dynamical systems generated by systems of linear random delay differential equations of the form
\begin{equation}\label{main-delay-eq}
z'(t) = A(\theta_{t}\w) \,z(t) + B(\theta_{t}\w) \,z(t-1),
\end{equation}
where $z(t) \in \R^N$, $N \ge 2$,  $A(\w)$, $B(\w)$ are  $N\times N$ real matrices:
\begin{equation*}
A(\w) =
\left(\begin{smallmatrix}
a_{11}(\w)&a_{12}(\w)&\cdots&a_{1N}(\w) \\
a_{21}(\w)&a_{22}(\w)&\cdots&a_{2N}(\w) \\
\vdots & \vdots & \ddots & \vdots \\
a_{N1}(\w)&a_{N2}(\w)&\cdots&a_{NN}(\w)
\end{smallmatrix}\right)\,,
\quad
B(\w) =
\left(\begin{smallmatrix}
b_{11}(\w)&b_{12}(\w)&\cdots&b_{1N}(\w) \\
b_{21}(\w)&b_{22}(\w)&\cdots&b_{2N}(\w) \\
\vdots & \vdots & \ddots & \vdots \\
b_{N1}(\w)&b_{N2}(\w)&\cdots&b_{NN}(\w)
\end{smallmatrix}\right)\,,
\end{equation*}
and $(\OFP,(\theta_t)_{t \in \R})$ is an ergodic metric dynamical system, with $\PP$ complete.\par\smallskip
From now on, the Euclidean norm on $\R^N$ will be denoted by $\norm{\cdot}$, $\R^{N \times N}$ will stand for the algebra of $N \times N$ real matrices with the operator or matricial norm induced by the Euclidean norm, i.e., $\norm{A} := \sup\{\norm{A\,u}: \norm{u} \le 1\}$, for any $A \in \R^{N \times N}$.  For $z = (z_1, \ldots, z_N) \in \R^N$ we write $z \gg 0$ if $z_i > 0$ for all $i \in \{1, 2, \ldots, N\}$.
\par\smallskip
For $1<p<\infty$, let $L=\R^N\times L_p([-1,0],\R^N)$ be the separable Banach space with the norm
\begin{equation*}
\normL{u}=\norm{u_1}+\norm{u_2}_p=\norm{u_1}+\left(\int_{-1}^0 \norm{u_2(s)}^p\,ds\right)^{\!\!1/p}
\end{equation*}
for any $u=(u_1,u_2)$ with $u_1\in\R^N$ and $u_2\in L_p([-1,0],\R)$.
The positive cone \[L^+=\left\{u=(u_1,u_2)\in L: u_1 \geq 0 \text{ and } u_2(s)\geq 0 \text{ for Lebesgue-a.e.~} s\in[0,1]\right\}\] is normal and reproducing, and the dual  $L^*=\R^N\times L_q([-1,0],\R)$ with $1/q+1/p=1$ is also separable.
\par\smallskip
We denote by $C$ the Banach space $C([-1, 0], \R^N)$ of continuous $\R^N$-valued functions defined on $[-1,0]$, with the supremum norm (denoted by $\normC{\cdot}$).  The positive cone
\begin{equation*}
    C^{+} = \left\{\, u \in C : u(s) \geq 0 \text{ for all } s \in [0,1] \,\right \}
\end{equation*}
is normal and reproducing.
\par
Further, by $J$ we denote the linear mapping from $C$ to $L$
\begin{equation}\label{defiJ}
\begin{array}{cccc}
J\colon & C & \longrightarrow & L\\
&  u &  \mapsto & (u(0),u) \,,
\end{array}
\end{equation}
which belongs to $\mathcal{L}(C,L)$ and $\norm{J}=2$.
\par\smallskip
Now we introduce the assumptions on the coefficients of the family~\eqref{main-delay-eq}:
\begin{enumerate}[label=(\textbf{S\arabic*}),series=supo,leftmargin=24pt]
\item\label{S1} (Measurability) $A, B \colon \Omega \to \R^{N \times N}$ are $(\mathfrak{F}, \mathfrak{B}(\R^{N \times N}))$\nobreakdash-\hspace{0pt}measurable.\smallskip
\item\label{S2} (Summability) The $(\mathfrak{F}, \mathfrak{B}(\R))$-measurable functions $a$, $b\colon \Omega \to \R$ defined as $a(\w):=\norm{A(\w)}$ and $b(\w):=\norm{B(\w)}$ have the properties:
\begin{align*}
	& \bigl[ \, \Omega \ni \w \mapsto a(\w) \in \R \, \bigr] \in L_1\OFP, \text{ and}
	\\
	& \Bigl[ \, \Omega \ni \w \mapsto \lnplus{\int_{0}^{1}b^q(\theta_{s}\w) \, ds} \in \R \, \Bigr] \in L_1\OFP.
\end{align*}
\end{enumerate}
As shown in~\cite[Remark 4.1]{MiNoOb2}
the following is sufficient for the fulfillment of the second condition in~\ref{S2}:
\begin{equation*}
\bigl[ \, \Omega \ni \w \mapsto b(\w) \in \R \, \bigr] \in L_q \OFP.
\end{equation*}
\begin{nota}\label{remark:a,b_L_loc}
Although the coefficients $A$ and  $B$ are defined only $\PP$-a.e. on $\Omega$, we can put the value of $A(\w)$ and $B(\w)$ to be equal to 0  for $\w$ in a set of null measure (see~\cite{MiNoOb2} for more details) to obtain
\begin{align*}
\bigl[\, \R\ni t\mapsto a(\theta_t\w)\in \R\,\bigr]&\in L_{1,\text{loc}}(\R)\,,\\
\bigl[ \, \R\ni t\mapsto b(\theta_t\w)\in \R\,\bigr]& \in L_{q,\text{loc}}(\R)\subset L_{1,\text{loc}}(\R) \,,
\end{align*}
for each $\w \in \Omega$.
\end{nota}
As a consequence, for a fixed $\w\in\Omega$ we will denote by  $U_{\w}^{0}(\cdot)$  the fundamental matrix solution of the system of Carath\'{e}odory linear ordinary  differential equations $z' = A(\theta_{t}\w)\,z$ and define
\begin{equation}\label{deficd}
c(\w) := \sup\limits_{0 \le t_1 \le t_2 \le 1} \norm{U^{0}_{\theta_{t_1}\w}(t_2 - t_1)}\,,\quad d(\w):=\biggl(\int_{-1}^0 \!b^q(\theta_{s+1}\w)\,ds\biggr)^{\!\!1/q}\,.
\end{equation}
Notice that $c(\w)\geq 1$, and as shown in~\cite[Lemma 4.2]{MiNoOb2} we have
\begin{equation}\label{ineqc}
c(\w) \le \exp\biggl( \int_{0}^{1} a(\theta_{s}\w) \, ds  \biggr).
\end{equation}
Under assumptions~\ref{S1} and \ref{S2}, it is shown in~\cite{MiNoOb2} that the family of systems~\eqref{main-delay-eq} generates \mlsps s both on $C$ and $L$.  More precisely, for $C$ we consider the initial value problem of Carath\'{e}odory type
\begin{equation}
\label{eq:IVP-C}
\begin{cases}
z'(t) = A(\theta_{t}\w) \,z(t) + B(\theta_{t}\w) \,z(t-1), & t \in [0, \infty)
\\
z(t) = u(t), & t \in [-1, 0],
\end{cases}
\end{equation}
where the initial datum $u$ is assumed to belong to $C=C([-1,0],\R^N)$, and to emphasize the dependence of the equation (resp.\ the initial value problem) on $\w \in \Omega$ we will write \eqref{main-delay-eq}$_{\w}$ (resp.\ \eqref{eq:IVP-C}$_{\w}$).  As shown in~\cite{MiNoOb2}, it has a unique solution denoted by $z(t,\w,u)$.  Moreover, it can be checked that for each $t$ and $r\geq 0$ \begin{equation}\label{cocycle}
z(t+r,\w,u)=z(t,\theta_r\w,z_r(\w,u))\,,
\end{equation}
where $z_r(\w,u)\colon [-1,0]\to \R$ is defined by $s\mapsto z(r+s,\w,u)$, and $z_t(\w,u)\in C$ for each $t\geq 0$, $\w\in\Omega$ and $u\in C$.
 Therefore, we can define the linear~operator
\begin{equation}\label{defiskpd}
\begin{array}{lccc}
 U^{(C)}_\w(t)\colon & C &\longrightarrow & C \\[.1cm]
                        & u & \mapsto & z_t(\w,u)\,.
\end{array}
\end{equation}
\par
Analogously, concerning $L$ we consider the initial value problem
\begin{equation}\label{eq:IVP}
\begin{cases}
z'(t) = A(\theta_{t}\w) \,z(t) + B(\theta_{t}\w)\, z(t-1), & t \in [0, \infty)
\\
z(t) = u_2(t), & t \in [-1, 0),
\\
z(0) = u_1,
\end{cases}
\end{equation}
with initial datum $u=(u_1,u_2)$ belonging to $L=\R^N\times L_p([-1,0],\R^N)$. Again, to emphasize the dependence of the initial value problem on $\w \in \Omega$ we will write \eqref{eq:IVP}$_{\w}$. It has a unique solution denoted by $z(t,\w,u)$ and, as in $C$ we can define the linear operator
\begin{equation*}
\begin{array}{lccc}
 U^{(L)}_\w(t)\colon & L &\longrightarrow & L \\[.1cm]
                        & u & \mapsto & (z(t,\w,u),z_t(\w,u))\,.
\end{array}
\end{equation*}
\par
Finally, under assumptions \ref{S1} and \ref{S2}, $\bigl((U^{(C)}_{\w}(t))_{\w \in \Omega, t \in \R^+)}, (\theta_{t})_{t \in \R}\bigr)$ and  $\bigl((U^{(L)}_{\w}(t))_{\w \in \Omega, t \in [0,\infty)}, (\theta_{t})_{t \in \R}\bigr)$ are  measurable linear skew\nobreakdash-\hspace{0pt}product semiflows on  $C$ and $L$ covering $\theta$, as shown in~\cite[Proposition 4.4 and 4.14, resp.]{MiNoOb2}. Moreover,  we can connect both semidynamical systems in the following way (see~\cite[Proposition~4.7]{MiNoOb2}). For $t\geq 1$ the linear operator
\begin{equation}\label{5.defiskpdLC}
\begin{array}{lccc}
 U^{(L,C)}_\w(t)\colon & L &\longrightarrow & C \\[.1cm]
                        & u & \mapsto & z_t(\w,u)
\end{array}
\end{equation}
belongs to $ \mathcal{L}(L, C)$ and  is a compact operator satisfying
\begin{equation}\label{eq:L-C-L}
	U^{(L,C)}_\w(t)= U_{\theta_{1}\w}^{(C)}(t-1)\circ U_{\w}^{(L,C)}(1)\,.
\end{equation}
Moreover, from~\cite[Corollary 4.8]{MiNoOb2} we know that
\begin{equation}\label{relationUC-UL-ULC}
U_\w^{(L)}(t)=J\circ U_\w^{(L,C)}(t)\,, \quad U_\w^{(C)}(t)=U_\w^{(L,C)}(t)\circ J
\end{equation}
for any $t\geq 1$ and any $\w\in\Omega$.
\par\smallskip
For the rest of the section the following assumptions will also be necessary.
\begin{enumerate}[resume*=supo]
\item\label{S3} (Cooperativity)\vspace{3pt}
\begin{itemize}
\item[(i)] $a_{ij}(\w)\geq 0$ for all $i\neq j$, $i,j=1,2,\ldots, N$ and $\w\in\Omega$.\vspace{3pt}
\item[(ii)] $b_{ij}(\w)\geq 0$ for all  $i,j=1,2,\ldots, N$ and $\w\in\Omega$.\vspace{3pt}
\end{itemize}
\item\label{S4} (Irreducibility) There is a $(\mathfrak{F}, \mathfrak{B}(\R))$-measurable function $\delta \colon \Omega \to (0,1]$, an invariant set $\widetilde\Omega_0$ of full measure $\PP(\widetilde\Omega_0)=1$  and  $M\in \N$ such that \smallskip
    \begin{itemize}
    \item[(i)] for each $\w\in\widetilde\Omega_0$ and $i\in \{1,2,\ldots,N\}$ there is a path  $j_1, j_2,j_3,\ldots,j_N$ starting  at $j_1=i$ such that $\{j_1,j_2,\ldots,j_N\}=\{1,2,\ldots, N\}$ and
    \[ \hspace{1.3cm} \int_0^M (a_{j_{l+1},j_l}(\theta_{s+t}\w)+b_{j_{l+1},j_l}(\theta_{s+t}\w))\,ds \geq \delta(\theta_t\w) \]
      whenever  $t\in\{k+(k-2)M:k=2,3,\dots,N\}$ and  $l=1,2,\ldots, N-1$.
    \item[(ii)] $\lnplus \ln (1/\delta)\in L_1\OFP$. \smallskip
\end{itemize}
\item\label{S5} \ref{S4}(i) holds and condition (ii)  is changed to \vspace{3pt}
	\begin{itemize}
	\item[(ii)] $\ln (1/\delta)\in L_1\OFP$.
	\end{itemize}
\end{enumerate}
\par\smallskip
Notice that~\ref{S5} implies~\ref{S4}. Under assumptions \ref{S1}--\ref{S4}  we will prove the existence  of a family of generalized principal Floquet subspaces   for
	the \mlsps\ $\bigl((U^{(L)}_{\w}(t))_{\w \in \Omega, t \in \R^+}, (\theta_{t})_{t \in \R}\bigr)$ generated by~\eqref{main-delay-eq}.  If, in addition,  \ref{S5} holds, then a generalized exponential separation of type II is obtained. In order to do that we first check the following results.  We call a solution $z(t)$ of~\eqref{main-delay-eq}$_\w$ \emph{positive} if, for any $t$ in its domain, $z_i(t) \ge 0$ for all $i \in  \{1,2,\ldots N\}$ and $z(t) \ne 0$.
\begin{prop}\label{TimeT}  Assume \textup{\ref{S1}}--\textup{\ref{S4}} and let $z(t)$ be a positive solution of~\eqref{main-delay-eq}$_\w$ satisfying $\,\sup_{t\in[0,1]} \norm{z(t)}>0$. Then   $z(t)\gg 0$ for each $t\geq N+(N-1)M$.
\end{prop}
\begin{proof}
Let $t^*$ be the point of $[0,1]$ and $j_1$ be the index in which the maximum is attained,  that is, $c_1=z_{j_1}(t^*)\geq  z_i(t)$ for each $i\in \{1,2,\ldots N\}$ and $t\in[0,1]$.  Consider the path and the constant $M$ given in \ref{S4}(i), and denote
\begin{equation}\label{Kj}
K_i(\w)= \exp\bigg(-\int_0^{N+H}\!\!\!|a_{j_ij_i}(\theta_s\w)|\,ds \bigg),\; \text{ for } i=1,2,\ldots,N
\end{equation}
where  $H=(N-1)M+1$.
From $z'_{j_1}(t)\geq a_{j_1j_1} (\theta_{t}\w) \,z_{j_1}(t)$, and using a comparison result for Carath\'{e}odory differential equations~\cite[Theorem 2]{paper:Walter}, we deduce that
\[z_{j_1}(t)\geq \exp\bigg(\int_{t^*}^{t}\!\!\!a_{j_ij_i}(\theta_s\w)\,ds\bigg)z_{j_1}(t^*)\geq c_1K_1(\w) \;\text{ for }  t\in[1,N+H]\,.\]
Analogously, if  $t\in[2+M,N+H]$,  from
\begin{align*}
z_{j_2}'(t) &\geq a_{j_2j_2}(\theta_t\w)\,z_{j_2}(t)+ a_{j_2j_1 }(\theta_t\w) z_{j_1}(t)+b_{j_2j_1}(\theta_t\w) z_{j_1}(t-1) \\
& \geq a_{j_2j_2}(\theta_t\w)\,z_{j_2}(t)+ \left[a_{j_2j_1 }(\theta_t\w) +b_{j_2j_1}(\theta_t\w) \right] c_1 K_1(\w)\,
\end{align*}
we obtain
\begin{align*}
z_{j_2}(t)&\geq  c_1K_1(\w)K_2(\w)\int_2^{t} \left[a_{j_2j_1 } +b_{j_2j_1}\right](\theta_{s}\w) \,ds\\
 &\geq c_1 K_1(\w)K_2(\w)\int_0^{M}\left[a_{j_2j_1 } +b_{j_2j_1}\right](\theta_{s+2}\w) \,ds
\end{align*}
 and \ref{S4}(i) provides
\[ z_{j_2}(t)\geq c_1 K_1(\w) K_2(\w)\,\delta(\theta_2\w)\;\text{ for }  t\in[2+M,N+H]\,.\]
Similarly,
\[ z_{j_3}(t)\geq c_1 K_1(\w) K_2(\w)K_3(\w)\,\delta(\theta_2\w)\,\delta(\theta_{3+M}\w)\;\text{ for }  t\in[3+2\,M,N+H]\,,\]
and finally, in a recursive way we prove  that for $k=1,2,\ldots N$
\begin{equation}\label{lowerIneq}
z_{j_k}(t)\geq c_1 \prod_{j=1}^{k} K_j(\w)\prod_{j=2}^{k}\delta(\theta_{j+(j-2)M} \w)\ \;\text{ for }  t\in[k+(k-1)M,N+H]\,,
\end{equation}
which finishes the proof.
\end{proof}
\begin{teor}\label{A1A3}
	Consider  the measurable linear skew-product semidynamical system $\bigl((U^{(L)}_{\w}(t))_{\w \in \Omega, t \in \R^+}, (\theta_{t})_{t \in \R}\bigr)$.  Under assumptions \textup{\ref{S1}}--\textup{\ref{S4}}, conditions \textup{\ref{A1}} and \textup{\ref{A3}}  hold for  time $T=N+(N-1)M+1$.
\end{teor}
\begin{proof}
	 As in~\cite[Proposition 5.7]{MiNoOb1}, concerning the first part of \ref{A1} with $T$  note that
	\[\sup\limits_{0 \le s \le T} {\lnplus{\n{U^{(L)}_{\w}(s)}}}\leq \sum_{k=0}^{T-1}\sup\limits_{k\le s \le k+1} {\lnplus{\n{U^{(L)}_{\w}(s)}}}\]
	and from cocycle property~\eqref{eq-cocycle}
	\begin{align*}
	\sup\limits_{k\le s \le k+1} \lnplus{\n{U^{(L)}_{\w}(s)}}&=\sup\limits_{0 \le s \le 1} {\lnplus{\n{U^{(L)}_{\w}(k+s)}}} \\
	&\leq \sup\limits_{0 \le s \le 1} \bigl(\lnplus\n{U^{(L)}_{\theta_s\w}(k)}+ \lnplus \n{U^{(L)}_{\w}(s)}\bigr).
	\end{align*}
	Therefore,
	\[\sup\limits_{0 \le s \le T} {\lnplus{\n{U^{(L)}_{\w}(s)}}}\leq T\,\sup\limits_{0 \le s \le 1} {\lnplus{\n{U^{(L)}_{\w}(s)}}}+ \sum_{k=1}^{T-1}\sup\limits_{0 \le s \le 1} {\lnplus{\n{U^{(L)}_{\theta_s\w}(k)}}}\]
	and we have to check that both terms belong to $L_1\OFP$.
As shown in~\cite[Proposition 4.6]{MiNoOb2},
$\norm{U^{(L)}_\w(t)}\leq 3\,c(\w)\,(1+d(\w))$   for  each $t\in[0,1]$  and  $\w\in\Omega$,
where $c(\w)$ and $d(\w)$ are defined in~\eqref{deficd} and $c(\w)\geq 1$.  Thus,  from   inequality~\eqref{ineqc} and \cite[Lemma 5.6]{MiNoOb1} we deduce that
\[\sup\limits_{0 \le s \le 1}{\lnplus{\n{U^{(L)}_{\w}(s)}}}\leq \ln 6 + \lnplus  c(\w)+ \lnplus d(\w) \leq \ln 6 +  \int_{0}^{1} a(\theta_{s}\w) \, ds +  \lnplus  d(\w)\,,\]
which belongs to $L_1\OFP$ because of the definition of $d$, \ref{S1}, \ref{S2}, Fubini's theorem and the invariance of $\PP$.
We omit the proof of the second term, as well as that of the second part of~\ref{A1} because they are analogous.  In addition, it is also immediate to check that~\ref{A2} follows from~\ref{S3}.
\par\smallskip
We will finish by verifying that \ref{A3} holds for time $T$.  First we claim that  for each  $\w\in\Omega$ and $u \in C^+$  with $U_\w^{(C)}(T)\, u\neq 0$
\begin{equation}\label{focusingC}
\widetilde\beta(\w,u)\,\widetilde\be\leq U_\w^{(C)}(T)\,u\leq \varkappa(\w)\,\widetilde \beta(\w,u)\,\widetilde\be
\end{equation}
where  $\widetilde \be\in C^+$ is the constant unit function $\widetilde\be(s)=1$ for each $s\in[-1,0]$,
\begin{align}
\widetilde\beta(\w,u)&=\normC{U_\w^{(C)}(1)\,u}\,k_\delta(\w)\,,\nonumber\\
\varkappa(\w)&=k_\delta^{-1}(\w)\,3^{T-1} \displaystyle\prod_{j=0}^{T-2}c(\theta_{j+1}\w)(1+d(\theta_{j+1}\w))\,, \label{varkappa}\\
k_\delta(\w)&=\prod_{j=1}^{N} K_j(\w)\prod_{j=2}^{N}\delta(\theta_{j+(j-2)M} \w)\,,\label{productKjdelta}
\end{align}
and $K_j(\w)$ are defined in~\eqref{Kj} ($N+H=T$).  Notice that $\normC{U_\w^{(C)}(1)\,u}$  is strictly positive  because of  $U_\w^{(C)}(T)\, u\neq 0$ and the cocycle property~\eqref{eq-cocycle}, so that $\widetilde\beta(\w,u)$ is strictly positive. Moreover, since $(U_\w^{(C)}(T)\,u)(s)= z(T+s,\w,u)$ for each $s\in[-1,0]$, the lower inequality of~\eqref{focusingC} follows from~\eqref{lowerIneq} with  $k=N$ and $t=T$ and \eqref{productKjdelta}.\par\smallskip
Concerning the upper inequality,  we first claim that for each $n\geq 1$
\begin{equation}\label{upper}
\normC{U_\w^{(C)}(t)\,u}\leq 3^n  \displaystyle\prod_{j=0}^{n-1}c(\theta_{j}\w)(1+d(\theta_{j}\w)) \,\normC{u}\;\text{ whenever }t\in[n-1,n]\,.\end{equation}
From~\cite[Proposition 4.4]{MiNoOb2}, we deduce that
\[
\normC{U^{(C)}_\w(t)\,u}\leq 3\,c(\w)\,(1+d(\w))\,\normC{u}  \quad\text{ for  each } t\in[0,1] \text{ and } \w\in\Omega\,,
\]
which together with the cocycle property~\eqref{eq-cocycle} yields
\begin{align*}
\normC{U_\w^{(C)}(t)\,u}&=\normC{U_{\theta_1\w}^{(C)}(t-1)\big(U_\w^{(C)}(1)\,u\big)}\\
 &\leq 3\,c(\theta_1\w)(1+d(\theta_1\w)\,3\,c(\w)(1+d(\w))\,\normC{u}
\end{align*}
for each $t\in[1,2]$ and $\w\in\Omega$. Thus~\eqref{upper} is easily checked in a recursive way.\par\smallskip
Finally, from  $U_\w^{(C)}(T)\,u\leq \normC{U_\w^{(C)}(T)\,u}\, \widetilde \be $  and~\eqref{upper}  for $t=T-1$ we obtain
\begin{align*}
U_\w^{(C)}(T)\,u&\leq\normC{U_{\theta_1\w}^{(C)}(T-1) (U_\w^{(C)}(1)\,u)}  \,\widetilde \be \\
& \leq \,\normC{U_\w^{(C)}(1)\,u)}\,3^{T-1} \displaystyle\prod_{j=0}^{T-2}c(\theta_{j+1}\w)(1+d(\theta_{j+1}\w))\,\widetilde \be\\
&=\varkappa(\w)\,\widetilde\beta(\w,u)\,\widetilde\be\,,
\end{align*}
as claimed.\par\smallskip
Next we will prove that for any $u\in L^+$ such that $U_\w^{(L)}(T)\,u\not=0$
\begin{equation}\label{focusingL}
\beta(\w,u)\,\be\leq U_\w^{(L)}(T)\,u\leq \varkappa(\w)\, \beta(\w,u)\,\be\,,
\end{equation}
where  $\be=(1/(2\sqrt{N}))\,(\widetilde\be(0),\widetilde\be)\in L^+$ is a unitary vector of $L$, i.e   $\normL{\be}=1$ and
\begin{equation}\label{betaL}
\beta(\w,u)=2\sqrt{N} \,\normC{U_\w^{(L,C)}(1)\,u}\, k_\delta(\w) >0\,.
\end{equation}
\par
As shown in~\cite{MiNoOb2}, the following relations hold
\begin{equation*}
\label{eq:L-C-2}
U^{(L)}_{\w}(t) = J \circ U^{(L, C)}_{\w}(t), \quad U^{(L, C)}_{\w}(t)=U_{\theta_1\w}^{(C)}(t-1)\circ U_\w^{(L,C)}(1)
\end{equation*}
for any $t \ge 1$ and any $\w \in \Omega$, where $J$ is the linear map defined in~\eqref{defiJ}. Hence
\[
U_\w^{(L)}(T)\,u=(z(T,\w,u),z_T(\w,u))=J\big( U_{\theta_1\w}^{(C)}(T-1)(U_\w^{(L,C)}(1)\,u)\big)
\]
and  notice that $U_\w^{(L,C)}(1)\,u=z_1(\w,u)\in C$. Therefore, as in Proposition~\ref{TimeT},
\[
z(T+s,\w,u)\geq \normC{U_\w^{(L,C)}(1)\,u} \,k_\delta(\w)\, \widetilde\be
\]
for each $s\in[-1,0]$ and the lower inequality of~\eqref{focusingL} holds.
Concerning the upper inequality, as above
\[z_T(\w,u)\leq \,\normC{U_\w^{(L,C)}(1)\,u)}\,3^{T-1} \displaystyle\prod_{j=0}^{T-2}c(\theta_{j+1}\w)(1+d(\theta_{j+1}\w))\,\widetilde \be\]
from which the upper inequality of~\eqref{focusingL} can be easily checked. \par\smallskip
In order to finish the proof we have to check that the $(\mathfrak{F},\mathfrak{B}(\R))$\nbd-measurable function $\varkappa$ defined by~\eqref{varkappa} satisfies
 $\lnplus\ln\varkappa\in L_1\OFP$.
From~\eqref{Kj} and the definition of $a(\w)=\norm{A(\w)}$ we deduce that there is a constant $c_0$ such that
\begin{equation}\label{ineqKj}
\Bigg(\prod_{j=1}^{N} K_j(\w))\Bigg)^{\!\!-1}\!\!\!= \exp\int_0^T \sum_{j=1}^N|a_{jj}(\theta_s\w)|\,ds\leq \exp\int_0^T c_0\,a(\theta_s\w) \,ds\,.
\end{equation}
  Moreover, from \eqref{ineqKj}, \eqref{ineqc} and $\ln(1+x)\leq x$ if $x\geq 0$
\[\lnplus \ln c(\theta_{j+1}\w)\leq \int_0^1 a(\theta_{j+1+s}\w)\,ds\, ,\quad \lnplus \ln(1+d(\theta_{j+1}\w))\leq \lnplus d(\theta_{j+1}\w)\,,\]
\begin{equation*}
0\leq \lnplus \ln\Bigg(\prod_{j=1}^{N} K_j(\w))\Bigg)^{\!\!-1}\!\!\!\leq c_0 \int_0^T a(\theta_s\w) \,ds\,,
\end{equation*}
and we deduce from~\cite[Lemma 5.6]{MiNoOb1} that there is an integer number $n_0$ such that
\begin{align*}
\lnplus \ln\varkappa (\w)\leq & \; c_0 \int_0^T a(\theta_s\w)\,ds+ \sum_{j=0}^{T-2} \left[\int_0^1 a(\theta_{j+1+s} \w)\,ds+ \ln ^+d(\theta_{j+1}\w)\,\right] \\[-.2cm]
 &\quad + \sum_{j=2}^N\lnplus \ln \left(1/\delta\right)(\theta_{j+(j-2)M}\w)+ \lnplus \ln(3^{T-1})\,+\ln (n_0).
\end{align*}
As before, from the definition of $d$, \ref{S1}, \ref{S2},  \ref{S4}(ii), Fubini's theorem and the invariance of $\PP$ we conclude that  $\lnplus \ln \varkappa\in L_1\OFP$  and \ref{A3} holds for $T$, as stated.
\end{proof}
 As a consequence, from~Theorem~\ref{Theorem3.10} the existence of a family of generalized principal Floquet subspaces is obtained.
\begin{teor}\label{generalized}
Under assumptions \textup{\ref{S1}}--\textup{\ref{S4}}
there is a family of generalized principal Floquet subspaces for the \mlsps\ $\Phi ^{(L)}= \allowbreak ((U_\w^{(L)}(t))_{\w \in \Omega, t \in \R^{+}}, \allowbreak (\theta_t)_{t\in\R})$ generated by~\eqref{main-delay-eq}, with generalized principal Lyapunov exponent  $\widetilde{\lambda}^{(L)}_1$.
\end{teor}
If we assume \textup{\ref{S5}} instead of \textup{\ref{S4}}, the existence of a generalized exponential separation of type II is obtained, as claimed before.
\begin{teor}\label{expsepL}
Under assumptions \textup{\ref{S1}}--\textup{\ref{S3}}  and  \textup{\ref{S5}}, the \mlsps\ $\Phi ^{(L)}= \allowbreak ((U_\w^{(L)}(t))_{\w \in \Omega, t \in \R^{+}}, \allowbreak (\theta_t)_{t\in\R})$ generated by~\eqref{main-delay-eq} admits a generalized exponential separation of type \textup{II} with $\widetilde{\lambda}_1^{(L)}> - \infty$.	
\end{teor}
\begin{proof}
As $L=\R^N\times L_p([-1,0],\R^N)$  is an ordered separable Banach space with separable dual  $L^*=\R^N\times L_q([-1,0],\R^N)$  and positive cone  $L^+$  normal and reproducing,  the claim follows from Theorems~\ref{thm:exponential-separation}, \ref{A1A3} and~\ref{generalized}, once we check that  $\widetilde\lambda_1^{(L)}>-\infty$.\par\smallskip
Let $\be$ be the unitary vector of $L$ defined for the focusing inequality~\eqref{focusingL}.  From Proposition~\ref{TimeT} we deduce that $U_\w^{(L)}(T)\,\be\neq 0$ and hence $U_\w^{(L)}(T)\,\be \geq \beta(\w,\be)\,\be$ for each $\w\in\Omega$. Thus, from the cocycle property~\eqref{eq-cocycle}, \ref{A2}, and the monotonicity of the norm, it is easy to check that
\[
\normL{U_\w^{(L)}(n\,T)\,\be}\geq \prod_{k=0}^{n-1} \beta(\theta_{kT}\w,\be) \quad \text{ for each } n\in\N\,,
\]
and therefore
\[
\frac{\ln \normL{U_\w^{(L)}(n\,T)\,\be}}{nT}\geq \frac{1}{nT}\sum_{k=0}^{n-1} \ln \beta(\theta_{kT}\w,\be)\,.
\]
Finally,  note that the function $-\ln \beta(\cdot,\be)$ for $\beta$ defined in~\eqref{betaL} satisfies
\[
0\leq -\ln \beta(\w,\be)\leq -\sum_{j=1}^N \ln K_j(\w)+\sum_{j=2}^N\ln \left(1/\delta)(\theta_{j+(j-2)M}\w\right)+ c_1
\] for some constant $c_1$, and thus belongs to $L_1\OFP$ from~\eqref{ineqKj} and~\ref{S5}.
Then, an application of Birkhoff ergodic theorem
gives that for $\PP$-a.e. $\w\in\Omega$ there holds
\[
\lim_{n\to\infty}\frac{1}{n}\sum_{k=0}^{n-1} \ln \beta(\theta_{kT}\w,\be)=\int_\Omega \ln\beta(\theta_T\w',\be)\,d\PP(\w')> -\infty\,,
\]
hence
\[
\limsup_{t\to\infty} \frac{\ln \normL{U_\w^{(L)}(t)\,\be}}{t}>-\infty\,,
\]
and \eqref{principal-exponent} finishes the proof.
\end{proof}
\begin{lema}\label{omegaxL->C} Let  $U_\w^{(L,C)}(1)$ be the compact linear operator defined in~\eqref{5.defiskpdLC}. Then  the mapping
\[
[\, \Omega\times L\ni (\w,u) \mapsto  U_\w^{(L,C)}(1)\,u \in C \,] \text{ is }(\mathfrak{F}\otimes\mathfrak{B}(L), \mathfrak{B}(C))\text{\nobreakdash-\hspace{0pt}measurable}\,.
\]
\end{lema}
\begin{proof}
From~\cite[Lemma 4.51, pp. 153]{AliB} it is enough to check that it is a Carath\'{e}odory function, i.e. for each fixed $\w\in\Omega$, the map $L\to C$, $u\mapsto  U_\w^{(L,C)}(1)\,u=z_1(\w,u)$ is continuous, which is well known, and for each fixed $u\in L$, the map
\begin{equation}\label{measurablemap}
	[\, \Omega\ni \w  \mapsto  U_\w^{(L,C)}(1)\,u=z_1(\w,u)\in C\,]  \text{ is } (\mathfrak{F}, \mathfrak{B}(C))\text{\nobreakdash-\hspace{0pt}measurable}.
\end{equation}
In order to prove this, denoting by $u=(c,v)\in \R^N\times L_p([-1,0],\R^N)$ we can find a sequence of continuous functions $v_n\in C$ such that, $v_n(0)$ converges to $c$ in $\R^N$ and $v_n$ converges to $v$ in $L_p([-1,0],\R^N)$ as $n \uparrow\infty$. Therefore,  the measurability of the maps 	$[\, \Omega\ni \w  \mapsto  z_1(\w,v_n(0),v_n)=z_1(\w,v_n)\in C\,]$  (see~\cite[Lemma 4.11]{MiNoOb2})  for each $n\in\N$, the convergence of them to the map~\eqref{measurablemap} as  $n \uparrow\infty$
and~\cite[Corollary 4.29]{AliB} show the measurability and finishes the proof.
\end{proof}
 \begin{teor}\label{generalizedC}
Under assumptions \textup{\ref{S1}}--\textup{\ref{S4}}
 there is a family of generalized principal Floquet subspaces for the \mlsps\ $\Phi ^{(C)}= \allowbreak ((U_\w^{(C)}(t))_{\w \in \Omega, t \in \R^{+}}, \allowbreak (\theta_t)_{t\in\R})$ generated by~\eqref{main-delay-eq}, with generalized principal Lyapunov exponent  $\widetilde{\lambda}^{(C)}_1$ which coincide with $\widetilde{\lambda}^{(L)}_1$ of \textup{Theorem~\ref{generalized}}.
\end{teor}
\begin{proof}
	Theorem~\ref{generalized} states the existence of a family $\{E_1^{(L)}(\w)\}_{\w \in \widetilde\Omega_1}$ of generalized principal Floquet subspaces for $\Phi^{(L)}$.  For $\w \in \widetilde{\Omega}_1$, put $E_1^{(L)}(\w) := \spanned\{w^{(L)}(\w)\}$ with $\normL{w^{(L)}(\w)} = 1$.
	First notice that from~\eqref{5.defiskpdLC} the function
	\[
	U_{\theta_{-1}w}^{(L,C)}(1)w^{(L)}(\theta_{-1}\w)
	=(z_1(\theta_{-1}\w,w^{(L)}(\theta_{-1}\w))\in \ C\,,
	\]
	from~\eqref{relationUC-UL-ULC}    we deduce that
	\[
	U_{\theta_{-1}\w}^{(L)}(1)w^{(L)}(\theta_{-1}\w)=(z_1(\theta_{-1}\w,w^{(L)}(\theta_{-1}\w))(0), z_1(\theta_{-1}\w,w^{(L)}(\theta_{-1}\w))\in \R^N\times C\,,
	\]
	which  is proportional to $w^{(L)}(\w)$ because $U_{\theta_{-1}w}^{(L)}(1)E_1^{(L)}(\theta_{-1} w)=E_1^{(L)}(\w)$, and hence,  can be denoted by  $w^{(L)}(\w)=(w(\w)(0),w(\w))\in \R^N\times C$. In particular, $w(\w)$ is proportional to $U_{\theta_{-1}\w}^{(L,C)}(1)w^{(L)}(\theta_{-1}\w)\neq 0$.
	\par\smallskip
	Next,  we define    $w^{(C)}(\w) := U_{\theta_{-1}\w}^{(L,C)}(1) \, w^{(L)}(\theta_{-1}\w) / \normC{U_{\theta_{-1}\w}^{(L,C)}(1) \, w^{(L)}(\theta_{-1}\w)} $,  unitary and also proportional to $w(\w)$.
	We consider   the one-dimensional subspace of $C$ defined by $E_1^{(C)}(\w):=\spanned\{w^{(C)}(\w)\}$  and  we claim that
	$\{E_1^{(C)}(\w)\}_{\w \in \widetilde\Omega_1}$  is a family of generalized principal Floquet subspaces of  $\Phi^{(C)}$, i.e. conditions (i)--(iv) of  Definition~\ref{generalized-floquet-space} hold. First, it is easy to deduce from Lemma~\ref{omegaxL->C} that $w^{(C)}\colon\widetilde\Omega_1\to C^+\setminus \{0\}$ is  $(\mathfrak{F}, \mathfrak{B}(C))$\nobreakdash-\hspace{0pt}measurable and  thus, (i) holds.
	Condition (ii), that is, the invariance of $\{E_1^{(C)}(\w)\}_{\w \in \widetilde\Omega_1}$, follows from $E_1^{(C)}(\w) = \spanned\{w(\w)\}$, $w^{(L)}(\w)=(w(\w)(0),w(\w))$ and the invariance of $\{E_1^{(L)}(\w)\}_{\w \in \widetilde\Omega_1}$.  In order to check (iii), notice that
	from \cite[Proposition 5.2(2)]{MiNoOb2}  we deduce that
	\begin{equation} \label{lambda(L)}
		\widetilde\lambda_1^{(L)}=\lim_{t\to\infty} \frac{\ln \normL{U_\w^{(L)}(t)\,w^{(L)}(\w)}}{t}=\lim_{t\to\infty} \frac{\ln \normC{U_\w^{(L,C)}(t)\,w^{(L)}(\w)}}{t}\,.
	\end{equation}
	However,  as stated in~\eqref{relationUC-UL-ULC}, $U_\w^{(C)}(t)=U_\w^{(L,C)}(t)\circ J$  for $t\geq 1$,  then
	\begin{equation}\label{lambda(C)}
		\widetilde\lambda_1^{(L)}= \lim_{t\to\infty} \frac{\ln \normC{U_\w^{(C)}(t)\,w(\w)}}{t}=\lim_{t\to\infty} \frac{\ln \normC{U_\w^{(C)}(t)\,w^{(C)}(\w)}}{t}\,,
	\end{equation}
 	as needed,  and in particular the generalized principal Lyapunov exponents of $\Phi^{(L)}$ and $\Phi^{(C)}$ coincide. Finally,  for each $\w\in\widetilde\Omega_1$ with $U_\w^{(C)}(1)\,u\neq 0$ we know that $U_\w^{(L)}(1)\,Ju\neq 0$. Then, as in~\cite[Proposition 5.2(1)]{MiNoOb2}, it can be shown that
	\begin{equation}\label{lambdabL=lambdaC}
		\limsup_{t\to\infty} \frac{\ln \normC{U_\w^{(C)}(t)\,u}}{t}=\limsup_{t\to\infty} \frac{\ln \normL{U_\w^{(L)}(t)\,Ju}}{t}\leq \widetilde\lambda_1^{(L)}=\widetilde\lambda_1^{(C)}\,.
	\end{equation}
Therefore, (iv) holds and $\{E_1^{(C)}(\w)\}_{\w \in \widetilde\Omega_1}$  is a family of generalized principal Floquet subspaces of  $\Phi^{(C)}$, as claimed.
\end{proof}
Although the dual  space of $C=C([-1, 0], \R^N)$ is not a separable Banach space and thus, Theorem~\ref{thm:exponential-separation} does not apply, we will show the existence of  a generalized exponential separation of type \textup{II} with $\widetilde{\lambda}_1^{(C)} > - \infty$ for the \mlsps\ $\Phi^{(C)}=\bigl((U^{(C)}_{\w}(t))_{\w \in \Omega, t \in \R^+}, (\theta_{t})_{t \in \R}\bigr)$.
 \begin{teor}\label{expsepC}
 Under assumptions \textup{\ref{S1}--\ref{S3}}  and  \textup{\ref{S5}}, the \mlsps\ $\Phi^{(C)} = \allowbreak ((U_\w^C(t))_{\w \in \Omega, t \in \R^{+}}, \allowbreak (\theta_t)_{t\in\R})$ generated by~\eqref{main-delay-eq} admits a generalized exponential separation of type \textup{II} with $\widetilde{\lambda}_1^{(C)}> - \infty$.	
\end{teor}
\begin{proof}
In  order to prove that $\Phi^{(C)}$ admits a generalized exponential separation of type~II we need to check that the  family $\{E_1^{(C)}(\w)\}_{\w \in \widetilde\Omega_1}$ of Theorem~\ref{generalizedC}, defined by $E_1^{(C)}(\w)=\spanned\{w(\w)\}$,  satisfies conditions (i)--(iii) of Definition~\ref{def:exponential-separation}. \par
\smallskip
 First notice that from Theorems~\ref{generalized} and  \ref{thm:exponential-separation} we know, among other things, that there are an invariant set $\widetilde\Omega_1$ of full measure $\PP(\widetilde\Omega_1)=1$ and  a family of generalized principal Floquet subspaces
   $\{E_1^{(L)}(\w)\}_{\w \in \widetilde\Omega_1}$ with $E_1^{(L)}(\w) = \spanned\{w^{(L)}(\w)\}$ satisfying that  the family of projections $\{ P^{(L)}(\w)\}_{\w\in\widetilde\Omega_1}$, associated to the invariant decomposition $L=E_1^{(L)}(\w)\oplus F_1^{(L)}(\w)$, is strongly measurable and tempered, i.e
   \begin{equation}\label{tempL}
   \lim\limits_{t \to \pm\infty} \frac{\ln{\n{P^{(L)}(\theta_{t}\w)}}}{t} = 0 \qquad \PP\text{-a.e. on }\widetilde\Omega_1\,,
   \end{equation}
   and $F_1^{(L)}(\w)\cap L^+=\{u\in L^+ : U_\w^{(L)}(1)\,u=0\}$ for any $\w\in\widetilde\Omega_1$. We also know that  $\normL{w^{(L)}(\w)}=1$ for any $\w\in\widetilde\Omega_1$. As in Remark~\ref{rem:tempered-definition}, the family of projections onto $E_1^{(L)}(\w)$ will be denoted by $\{ \widetilde P^{(L)}(\w)\}_{\w\in\widetilde\Omega_1}$ (remember that $\widetilde{P}^{(L)}(\w) = \mathrm{Id}_L - P^{(L)}(\w)$). As $\widetilde{P}^{(L)}(\w)$ is a projection on $\spanned\{w^{(L)}(\w)\}$, for each $\widehat u\in L$ we have $ \widetilde P^{(L)}(\w)\,\widehat u=\lambda(\w,\widehat u)\, w^{(L)}(\w)$.
   \par\smallskip
Next, for any  $u\in C$, we know that $Ju=(u(0),u)\in L$ so that $u$ can be decomposed as $u= \lambda(\w,Ju)\,w(\w)+ \bigl( u-\lambda(\w,Ju)\,w(\w) \bigr)$, with $u-\lambda(\w,Ju)\,w(\w)\in C$.
For $\w \in \widetilde{\Omega}_1$, put $F_1^{(C)}(\w) := \{\, u \in C : J u \in F_1^{(L)}(\w) \,\}$. The closedness of $F_1^{(C)}(\w)$ in $C$ is a consequence of the continuity of $J$ and the closedness of $F_1^{(L)}(\w)$ in $L$.
\par\smallskip
We have thus obtained an invariant decomposition $C = E_1^{(C)}(\w) \oplus F_1^{(C)}(\w)$.  We claim that for any $\w \in \widetilde{\Omega}_1$ there holds $F_1^{(C)}(\w) \cap C^+ = \{u\in C^+ : U_\w^{(C)}(1)\,u=0 \}$.  Indeed, let $u \in C^{+}$ be such that $J u \in F_1^{(L)}(\w)$.  Then, by~\eqref{relationUC-UL-ULC} and the invariance of $F_1^{(L)}$, $J U_{\w}^{(C)}(1)\, u = J U_{\w}^{(L,C)}(1) J u = U_{\w}^{(L)}(1) J u \in F_1^{(L)}(\theta_{1}\w)$.  Further, it follows from~\ref{A2} and Definition~\ref{def:exponential-separation}(i)  for $\Phi^{(L)}$ that $J U_{\w}^{(C)}(1)\, u = 0$.  As $J$ is injective, $U_{\w}^{(C)}(1) \,u = 0$ holds.  On the other hand, if $u \in C^{+}$ is such that $U_{\w}^{(C)}(1)\, u = 0$, then, by~\eqref{relationUC-UL-ULC}, $0 = J U_{\w}^{(L,C)}(1) J u = U_{\w}^{(L)}(1) J u$, and, since $J u \in L^{+}$,  again by Definition~\ref{def:exponential-separation}(i)  for $\Phi^{(L)}$, $J u \in F_1^{(L)}(\w)$, that is, $u \in F_1^{(C)}(\w)$,  our claim is true and  (i) of Definition~\ref{def:exponential-separation} holds for $\Phi^{(C)}$.
\par\smallskip
Next, in view of Remark~\ref{rem:tempered-definition},  we prove that the family   $\{\widetilde P^{(C)}(\w)\}_{\w\in\widetilde\Omega_1}$   of projections onto $E_1^{(C)}(\w)$,  defined as $\widetilde{P}^{C}(\w)\,u=\lambda(\w,Ju)\,w(\w)$ for each $u\in C$,   is strongly measurable and tempered.  Concerning  the strong measurability, we have to check that for each $u \in C$ the mapping $[\, \widetilde \Omega_1 \ni \w \mapsto \widetilde P^{(C)}(\w)u \in C \,]$ is $(\mathfrak{F}, \mathfrak{B}(C))$\nobreakdash-\hspace{0pt}measurable. Once we  fix $u\in C$, the strong measurability of  $\{ \widetilde P^{(L)}(\w)\}_{\w\in\widetilde\Omega_1}$ together with  $ \widetilde P^{(L)}(\w)\, Ju=\lambda(\w,Ju)\, w^{(L)}(\w)$ and $\normL{w^{(L)}(\w)}=1$ show that $[\, \widetilde \Omega_1 \ni \w \mapsto \lambda(\w,Ju) \in \R \,]$ is $(\mathfrak{F}, \mathfrak{B}(\R))$\nobreakdash-\hspace{0pt}measurable, and the result follows from the composition of the maps
\[ \widetilde\Omega_1  \to  \R\times C  \to   C\,,\quad  \w \mapsto (\lambda(\w,Ju),w(\w))\mapsto \lambda(\w,Ju)\,w(\w)\,.\]
\par
In order to show that the family is tempered, first notice that from the definition of $\widetilde P^{(L)}$ and $\normL{w^{(L)}(\w)}=1$  we deduce that
\[
 \norm{\widetilde P^{(L)}(\theta_t\w)} = \sup_{\normL{\widehat u}\leq 1}\normL{ \widetilde P^{(L)}(\theta_t\w)\,\widehat u}= \sup_{\normL{\widehat  u}\leq 1} |\lambda(\theta_t\w,\widehat u)| \text{ for any }  \w\in\widetilde\Omega_1\,.
\]
Therefore, from  $\{Ju\in L : u\in C \text{ and } \normC{u}\leq 1\} \subset \{\widehat u\in L: \normL{\widehat u}\leq 2\}$  it holds
\begin{align*}
1\leq \norm{\widetilde P^{(C)}(\theta_t\w)} & =\sup_{\normC{u}\leq 1}\normC{ \widetilde P^{(C)}(\theta_t\w)\,u}= \sup_{\normC{u}\leq 1} |\lambda(\theta_t\w,Ju)\,|\normC{w(\theta_t\w)}\\
&\leq \sup_{\normL{\widehat u}\leq 2} |\lambda(\theta_t\w,\widehat u)| \,\normC{w(\theta_t\w)}\leq 2 \,\norm{\widetilde P^{(L)}(\theta_t\w)} \,\normC{w(\theta_t\w)}\,,
\end{align*}
and consequently,
\[
0\leq \lim\limits_{t \to \pm\infty} \frac{\ln{\n{\widetilde P^{(C)}(\theta_{t}\w)}}}{t}\leq  \lim\limits_{t \to \pm\infty} \frac{\ln{\n{\widetilde P^{(L)}(\theta_{t}\w)}}}{t}+  \lim\limits_{t \to \pm\infty} \frac{\ln{\normC{w(\theta_t\w)}}}{t}\,.
\]
Thus,  because of~\eqref{tempL}, in order to prove  that   $\{\widetilde P^{(C)}(\w)\}_{\w\in\widetilde\Omega_1}$   is tempered,  it is enough to check that
\begin{equation}\label{limit_w/t}
\lim\limits_{t \to \pm\infty} \frac{\ln{\normC{w(\theta_t\w)}}}{t}=0\,.
\end{equation}
From  $U_\w^{(L)}(t)\,w^{(L)}(\w)=\alpha(t,\w)\,w^{(L)}(\theta_t \w)$
and $\normL{w^{(L)}(\theta_tw)}=1$ we obtain
\[
w^{(L)}(\theta_t\w)=\frac{U_\w^{(L)}(t)\,w^{(L)}(\w)}{\normL{U_\w^{(L)}(t)\,w^{(L)}(\w)}} \, , \;\text{ and thus, }  w(\theta_t\w)=\frac{U_\w^{(L,C)}(t)\,w^{(L)}(\w)}{\normL{U_\w^{(L)}(t)\,w^{(L)}(\w)}}
\]
for $t\geq 0$. Hence,  $\ln{\normC{w(\theta_t\w)}}= \ln{\normC{U_\w^{(L,C)}(t)\,w^{(L)}(\w)}}-\ln \normL{U_\w^{(L)}(t)\,w^{(L)}(\w)}$ and~\eqref{limit_w/t} as $t$ goes to $+\infty$ follows from~\eqref{lambda(L)}.  Concerning the limit as $t$ goes to~$-\infty$,  as in Theorem~\ref{Theorem3.10} (see~\cite[Theorem 3.10]{MiNoOb1} for the complete proof)  we consider the negative semiorbit for $\Phi^{(L)}$ defined ~as
\begin{equation}
	\label{eq:definition-of-w}
	w_\w^{(L)}(s)
	=\frac{w^{(L)}(\theta_s\w)}{\normL{U_{\theta_s\w}^{(L)}(-s) \, w^{(L)}(\theta_s\w )}}  \quad \text{for } s\leq 0\,,
\end{equation}
 which satisfies $w_\w^{(L)}(0)= w^{(L)}(\w)$ and
\begin{equation}\label{lambda1L}
	  \widetilde\lambda_1^{(L)} = \lim_{s\to- \infty} \frac{1}{s}\ln\normL{w_\w^{(L)}(s)} =-\lim_{s\to-\infty} \frac{\ln \normL{U_{\theta_s \w}^{(L)}(-s)\,w^{(L)}(\theta_s\w)}}{s}\,.
\end{equation}
Then, as in~\cite[Theorem 5.6]{MiNoOb2},   denoting  by  $w_\w(s)$
the second component  of  $w_\w^{(L)}(s)$,  which belongs to $C$,  defines a negative semiorbit for $\Phi^{(C)}$,  $J w_\w(s)=  w_\w^{(L)}(s)$ and
\begin{equation}\label{limitws}
\lim_{s\to-\infty}\frac{\ln{\normC{w_\w(s)}}}{s}=\widetilde\lambda_1^{(L)}\,.
\end{equation}
Therefore, from~\eqref{eq:definition-of-w}  we deduce that
\[
w(\theta_s\w)=\normL{U_{\theta_s \w}^{(L)}(-s)\,w^{(L)}(\theta_s\w)}\,w_\w(s) \quad \text{ for each } s\leq 0\,,
\]
and   $\ln\normC{w(\theta_s\w)}= \ln\normL{U_{\theta_s \w}^{(L)}(-s)\,w^{(L)}(\theta_s\w)}+\ln \normC{w_\w(s)}$ holds. Thus from~\eqref{lambda1L} and \eqref{limitws}  the limit \eqref{limit_w/t}  holds as $t$ goes to $-\infty$,  $\{\widetilde P^{(C)}(\w)\}_{\w\in\widetilde\Omega_1}$   is tempered, as claimed, and (ii) of Definition~\ref{def:exponential-separation}~holds.\par
Now we check that (iii) of this definition also holds, that is,
there exists $\widetilde{\sigma}\in (0,\infty]$ such that
\[ \lim_{t\to\infty}\frac{1}{t}\ln \frac{\n{U_\w^{(C)}(t)|_{F_1^{(C)}(\w)}}}{\n{U_\w^{(C)}(t)\,w(\w)}_C}=-\widetilde{\sigma}
\quad  \text{for each } \w\in\widetilde\Omega_1\,.\]
Again, from Theorems~\ref{expsepL} and~\ref{thm:exponential-separation} we know that  there exists $\widetilde\sigma\in(0,\infty]$ such that  $\widetilde\lambda^{(L)}_2=\widetilde\lambda^{(L)}_1-\widetilde\sigma$ and
\[ \lim_{t\to\infty}\frac{1}{t}\ln \frac{\n{U_\w^{(L)}(t)|_{F_1^{(L)}(\w)}}}{\n{U_\w^{(L)}(t)\,w^{(L)}(\w)}_L}=-\widetilde{\sigma}
\quad \text{and} \quad \lim_{t\to\infty}\frac{1}{t}\ln\n{U_\w^{(L)}(t)|_{F_1^{(L)}(\w)}}= \widetilde\lambda^{(L)}_2
\]
for each $\w\in\widetilde\Omega_1$, which together with~\eqref{lambda(C)} shows that it is enough to prove that
\begin{equation}\label{lambda2C}
 \lim_{t\to\infty}\frac{1}{t}\ln\n{U_\w^{(C)}(t)|_{F_1^{(C)}(\w)}} = \widetilde\lambda^{(L)}_2  \quad  \text{for each } \w\in\widetilde\Omega_1\,.
 \end{equation}
As in \cite[Proposition 5.2(2)]{MiNoOb2}  it is easy to check that
\begin{equation}\label{lambda2otra}
 \widetilde\lambda^{(L)}_2 =\lim_{t\to\infty}\frac{1}{t}\ln\n{U_\w^{(L)}(t)|_{F_1^{(L)}(\w)}} = \lim_{t\to\infty}\frac{1}{t}\ln\n{U_\w^{(L,C)}(t)|_{F_1^{(L)}(\w)}}
\end{equation}
for each  $\w\in\widetilde\Omega_1$. Moreover, from~\eqref{eq:L-C-L} and the construction of $F^{(C)}_1(\w)$ we obtain
\begin{equation}\label{subsets}
U_\w^{(L,C)}(1)\big(F_1^{(L)}(\w)\big)\subset F_1^{(C)}(\theta_1\w) \quad\text{and}\quad J(F_1^{(C)}(\w))\subset F_1^{(L)}(\w)\,.
\end{equation}
 Thus,  from~\eqref{eq:L-C-L}, \eqref{relationUC-UL-ULC}, \eqref{lambda2otra} and  \eqref{subsets} we deduce the following chain of inequalities
\begin{align*} 
\widetilde\lambda^{(L)}_2&=\lim_{t\to\infty}\frac{1}{t+1}\ln\n{U_{\theta_{-1}\w}^{(L,C)}(t+1)|_{F_1^{(L)}(\theta_{-1}\w)}} \\
 &= \lim_{t\to\infty}\frac{1}{t+1}\ln\n{U_\w^{(C)}(t)\circ U_{\theta_{-1}\w}^{(L,C)}(1)|_{F_1^{(L)}(\theta_{-1}\w)}}\\
 & \leq \liminf_{t\to\infty}\frac{1}{t}\ln\n{U_\w^{(C)}(t)|_{F_1^{(C)}(\w)}}\leq \limsup_{t\to\infty}\frac{1}{t}\ln\n{U_\w^{(C)}(t)|_{F_1^{(C)}(\w)}} \\
 & = \limsup_{t\to\infty}\frac{1}{t}\ln\n{U_\w^{(L,C)}(t)\circ J|_{F_1^{(C)}(\w)}}
 \leq \lim_{t\to\infty}\frac{1}{t}\ln\n{U_\w^{(L,C)}(t)|_{F_1^{(L)}(\w)}} =  \widetilde\lambda^{(L)}_2,
\end{align*}
which shows~\eqref{lambda2C} and (iii) of Definition~\ref{def:exponential-separation}~holds. Finally,   from  $\widetilde{\lambda}_1^{(L)}> - \infty$  and ~\eqref{lambdabL=lambdaC} we deduce that $\widetilde{\lambda}_1^{(C)}> - \infty$, which finishes the proof.
\end{proof}
\section{Semiflows generated by linear random delay systems in a state space with not separable dual}\label{sect:lastsection}
This section deals with three cases in which the dual of the phase Banach space $X$ is non-separable and $\OFP$ is a Lebesgue space. In particular, we will apply the Oseledets theory of Section~\ref{sect:Oseledets}, more precisely Theorem~\ref{thm:Appendix}, to show the existence of a generalized exponential separation of type~II.
\subsection{Case 1} In this subsection we briefly explain what happens for $p=1$, that is, we consider the separable Banach space $\widehat{L}=\R^N\times L_1([-1,0],\R^N)$   with the norm
\begin{equation*}\label{normL+2a}
\norm{u}_{\widehat{L}}=\norm{u_1}+\norm{u_2}_{1}=\norm{u_1}+\int_{-1}^0 \norm{u_2(s)} \, ds
\end{equation*}
for any $u=(u_1,u_2)$ with $u_1\in \R^N$ and $u_2\in L_1([-1,0],\R^N)$.
We maintain conditions~\ref{S1}, \ref{S3}--\ref{S5} and we change  \ref{S2} by
\begin{enumerate}[label=(\textbf{S2\alph*}),series=supo,leftmargin=30pt]
\item\label{S2+a} The $(\mathfrak{F}, \mathfrak{B}(\R))$-measurable functions $a$, $b \colon \Omega \to \R$ defined as $a(\w) := \norm{A(\w)}$ and $b(\w):=\norm{B(\w)}$ have the properties:
\begin{align*}
& \bigl[ \, \Omega \ni \w \mapsto a(\w) \in \R \, \bigr] \in L_1\OFP, \text{ and}
\\
&
\bigl[ \, \Omega \ni \w \mapsto \lnplus  \esssup_{-1\le s\le 0} b(\theta_{s+1} \w) \in \R \, \bigr] \in L_{1} \OFP.
\end{align*}
\end{enumerate}
\begin{nota}\label{remark:a,b,L1}
Note that
\begin{equation*}
    \bigl[ \, \Omega \ni \w \mapsto  b(\w) \in \R \, \bigr] \in L_{\infty} \OFP
\end{equation*}
is sufficient for the second part of the assumption~\ref{S2+a}. Moreover, under assumption {\rm \ref{S2+a}}, we can deduce that  for $\PP$-a.e. $\w \in \Omega$
\begin{equation*}
\bigl[ \, \R\ni t\mapsto b(\theta_t\w)\in \R\,\bigr]  \in L_{\infty,\text{loc}}(\R).
\end{equation*}
\end{nota}
\par\smallskip
Under assumptions \ref{S1}, \ref{S2+a}, the initial value problem~\eqref{eq:IVP} with initial datum $u=(u_1,u_2)\in \widehat{L}$ has a unique solution $z(\cdot,\w,u)$, and as before for $L$, we will denote
$ U^{( \widehat{L})}_{\w}(t)\colon  \widehat{L} \longrightarrow  \widehat{L}$, $u  \mapsto  (z(t,\w,u),z_t(\w,u))$.
If we change $d(\w)$ in~\eqref{deficd} by
 \begin{equation}\label{widehatd}
 \widehat{d}(\w):= \esssup_{-1\le s\le 0} b(\theta_{s+1} \w)\,,
 \end{equation}
 it is not hard to check that Proposition~4.14 of~\cite{MiNoOb2} holds for the new fiber space
 $ \widehat{L}$, that is,  $\Phi^{(\widehat{L})}=\bigl((U^{(\widehat{L})}_{\w}(t))_{\w \in \Omega, t \in \RR^+}, (\theta_{t})_{t \in \R}\bigr)$ is a measurable linear skew-product semiflow on $\widehat{L}$ covering $\theta$. This result was already stated and used in~\cite[pp. 2254]{MiNoOb2}.
 \par\smallskip
Next  we prove the existence a family of generalized principal Floquet subspaces for this case.
\begin{teor}
Under assumptions \textup{\ref{S1}}, \textup{\ref{S2+a}}, \textup{\ref{S3}} and \textup{\ref{S4}} there is a family of generalized principal Floquet subspaces for the \mlsps\ $\Phi^{(\widehat{L})}=\bigl((U^{(\widehat{L})}_{\w}(t))_{\w \in \Omega, t \in \RR^+}, (\theta_{t})_{t \in \R}\bigr)$ generated by~\eqref{main-delay-eq}, with generalized principal Lyapunov exponent  $\widetilde{\lambda}^{(\widehat{L})}_1$.
\end{teor}
\begin{proof}
    First note that from~\cite[Lemma 5.8]{MiNoOb2} we have $\norm{U^{(\widehat{L})}_\w(t)}\leq 3\,c(\w)\,(1+\widehat d(\w))$ where $\widehat{d}$ is defined in~\eqref{widehatd}.
In addition, Proposition~\ref{TimeT} remains  true when~\ref{S2+a} is assumed instead of~\ref{S2}.  Hence, as in Theorem~\ref{A1A3} we  can deduce step by step  that  conditions~\ref{A1} and \ref{A3} hold for time $T=N+(N-1)M+1$  for the semiflow~$\Phi^{(\widehat{L})}$. Finally, Theorem~\ref{Theorem3.10} shows that  a family of generalized principal Floquet subspaces of $\Phi^{(\widehat L)}$ is obtained, as claimed.
\end{proof}
Finally, from Section~\ref{sect:Oseledets} an exponential separation of type \textup{II} for the \mlsps\ $\Phi^{(\widehat{L})}$ is obtained.
\begin{teor}
   Assume \textup{\ref{S1}}, \textup{\ref{S2+a}}, \textup{\ref{S3}} and  \textup{\ref{S5}}. Then the \mlsps\ $\Phi^{(\widehat{L})}=\bigl((U^{(\widehat{L})}_{\w}(t))_{\w \in \Omega, t \in \R^+}, (\theta_{t})_{t \in \R}\bigr)$ generated by~\eqref{main-delay-eq} admits a generalized exponential separation of type \textup{II} with $\widetilde{\lambda}_1^{(\widehat{L})}> - \infty$.
\end{teor}
\begin{proof}
From Theorem~\ref{thm:Appendix} it suffices to check that conditions~\ref{A1} and \ref{A3} hold, the operator $U^{(\widehat{L})}_{\w}(2)$ is compact and $\widetilde{\lambda}_1^{(\widehat{L})}> - \infty$. Since condition \textup{\ref{S4}} is weaker that \textup{\ref{S5}}, from the previous theorem we know that~\ref{A1} and \ref{A3} hold. Moreover,  the compactness of the operator $U^{(\widehat{L})}_{\w}(2)$ is shown in~\cite[Lemma 5.9(iii)] {MiNoOb2}. Finally, from Lemma~\ref{lm:top-equals-principal} we deduce that  $\widetilde{\lambda}_1^{(\widehat{L})}=\lambda_{\mathrm{top}}$, so it suffices to show that $\lambda_{\mathrm{top}}>-\infty$, which can be done via Birkhoff ergodic theorem as in Theorem~\ref{expsepL}, and finishes the proof.
\end{proof}
\subsection{Case 2} In this subsection we briefly explain what happens for the separable Banach space $C=C([-1,0],\RR^N)$ with the usual sup-norm when we maintain conditions~\ref{S1}, \ref{S3}--\ref{S5} and we replace \ref{S2} by
\begin{enumerate}[label=(\textbf{S2\alph*}),series=supo,leftmargin=30pt, resume]
\item\label{S2+b} The $(\mathfrak{F}, \mathfrak{B}(\R))$-measurable functions $a$, $b\colon \Omega \to \R$ defined as $a(\w):=\norm{A(\w)}$ and $b(\w):=\norm{B(\w)}$ have the properties:
\begin{align*}
& \bigl[ \, \Omega \ni \w \mapsto a(\w) \in \R \, \bigr] \in L_1\OFP, \text{ and}\\
& \Bigl[ \, \Omega \ni \w \mapsto \lnplus  \int_{-1}^0 \!b(\theta_{s+1}\w)\,ds \in \R \, \Bigr] \in L_{1} \OFP.
\end{align*}
\end{enumerate}
\begin{nota}\label{remark:a,b,C}
Note that
\begin{equation*}
     \bigl[ \, \Omega \ni \w \mapsto  b(\w) \in \R \, \bigr] \in L_{1} \OFP.
\end{equation*}
is sufficient for the second part of the assumption~\ref{S2+b}. Moreover, under assumption {\rm \ref{S2+b}}, we can deduce that  for $\PP$-a.e. $\w \in \Omega$
\begin{equation*}
\bigl[ \, \R\ni t\mapsto b(\theta_t\w)\in \R\,\bigr]  \in L_{1,\text{loc}}(\R).
\end{equation*}
\end{nota}
\par
As in Section~\ref{sect:semiflows}, it can be checked that, under assumptions \ref{S1}, \ref{S2+b},  for each $u\in C$ the initial value problem~\eqref{eq:IVP-C} has a unique Carath\'{e}odory type solution  which will be denoted by $z(\cdot,\w,u)$. Moreover, the cocycle relation~\eqref{cocycle} is also satisfied and we will use the same notation $U^{(C)}_\w(t)$ for the corresponding linear operator~\eqref{defiskpd}. Then,  a \mlsps\ is also obtained, as shows the next result.
\begin{prop}\label{prop:delay-semiflow-C}
Under \textup{\ref{S1}} and \textup{\ref{S2+b}}, $\bigl((U^{(C)}_{\w}(t))_{\w \in \Omega, t \in \mathbb{R}^{+}}, (\theta_{t})_{t \in \R}\bigr)$ is a measurable linear skew\nobreakdash-\hspace{0pt}product semiflow on $C$ covering $\theta$.
\end{prop}
\begin{proof}
First of all in~\eqref{deficd} we change $d$ by
\begin{equation*}
\widetilde{d}(\w):=\int_{-1}^0 \!b(\theta_{s+1}\w)\,ds.
\end{equation*}
From this it is not difficult to check that
\begin{equation}\label{UCbounded}
\normC{U^{(C)}_\w(t)\,u}\leq c(\w)\,(1+\widetilde{d}(\w))\,\normC{u}  \quad\text{ for  each } t\in[0,1] \text{ and } \w\in\Omega\,,
\end{equation}
and together with the cocycle property~\eqref{cocycle} we deduce that $U^{(C)}_\w(t)\in  \mathcal{L}(C)$ for each $t\geq 0$ and $\w\in\Omega$. The rest of the proof follows step by step the one of~\cite[Proposition 4.12]{MiNoOb2}.
\end{proof}
\begin{teor}\label{generalized-case2}
Under assumptions \textup{\ref{S1}}, \textup{\ref{S2+b}}, \textup{\ref{S3}} and \textup{\ref{S4}} there is a family of generalized principal Floquet subspaces for the \mlsps\ $\Phi^{(C)}=\bigl((U^{(C)}_{\w}(t))_{\w \in \Omega, t \in \RR^+}, (\theta_{t})_{t \in \R}\bigr)$ generated by~\eqref{main-delay-eq}, with generalized principal Lyapunov exponent  $\widetilde{\lambda}^{(C)}_1$.
\end{teor}
\begin{proof}
First note that Proposition~\ref{TimeT} remains  true when~\ref{S2+b} is assumed instead of~\ref{S2}. Hence, from~\eqref{UCbounded}, as in Theorem~\ref{A1A3}, we  can deduce that  conditions~\ref{A1} and \ref{A3} hold for time $T=N+(N-1)M+1$  for the semiflow~$\Phi^{(C)}$. Finally, as in Case 1, Theorem~\ref{Theorem3.10} finishes the proof.
\end{proof}
Finally, once we check that $U_\w^{(C)}(1)$  is a compact operator for each $\w\in\Omega$, from Section~\ref{sect:Oseledets} an exponential separation of type \textup{II} for the \mlsps\ $\Phi^{(C)}$ is obtained.
\begin{lema}\label{compactUCS1+S2+b}
 Under \textup{\ref{S1}} and \textup{\ref{S2+b}}, the bounded operator $U^{(C)}_{\w}(1)$ is compact for any $\w\in\Omega$.
\end{lema}
\begin{proof}
   First notice that the equicontinuity of the set $\{U^{(C)}_{\w}(1)\,u: \norm{u}_C\le 1\}$ follows from~\ref{S1}, \ref{S2+b}, \eqref{UCbounded}, and the inequality
\begin{align*}
 \norm{(U^{(C)}_{\w}(1)u)(s_1) & -(U^{(C)}_{\w}(1)u)(s_2)}_{C} \\ &\le \int_{1+s_1}^{1+s_2}  a(\theta_{s}\w)\norm{z(s,\w,u)}\,ds+\int_{s_1}^{s_2}  \
 b(\theta_{s+1}\w)\norm{u(s)}\,ds \\
&\leq c(\w)(1+\widetilde{d}(\w))\norm{u}_{C} \int_{1+s_1}^{1+s_2} a(\theta_{s}\w)\,ds+ \norm{u}_{C} \int_{s_1}^{s_2} b(\theta_{s+1}\w)\,ds.
\end{align*}
 whenever $-1\le s_1\le s_2 \le 0$.
Therefore from $U^{(C)}_{\w}(1)\in \mathcal{L}(C)$, the Ascoli--Arzel\`{a} theorem finishes the proof.
\end{proof}
\begin{teor}\label{expIIcase2}
   Assume \textup{\ref{S1}}, \textup{\ref{S2+b}}, \textup{\ref{S3}} and  \textup{\ref{S5}}. Then the \mlsps\ $\Phi^{(C)}=\bigl((U^{(C)}_{\w}(t))_{\w \in \Omega, t \in \R^+}, (\theta_{t})_{t \in \R}\bigr)$ generated by~\eqref{main-delay-eq} admits a generalized exponential separation of type \textup{II} with $\widetilde{\lambda}_1^{(C)}> - \infty$.
\end{teor}
\begin{proof}
From Theorem~\ref{thm:Appendix} it suffices to check that conditions~\ref{A1} and \ref{A3} hold, the operator $U^{(C)}_{\w}(1)$ is compact and $\widetilde{\lambda}_1^{(C)}> - \infty$. Since condition \textup{\ref{S4}} is weaker that \textup{\ref{S5}}, from the previous theorem and lemma we know that~\ref{A1} and \ref{A3} hold and $U^{(C)}_{\w}(1)$  is compact. Finally, from Lemma~\ref{lm:top-equals-principal} we deduce that  $\widetilde{\lambda}_1^{(C)}=\lambda_{\mathrm{top}}$, so it suffices to show that $\lambda_{\mathrm{top}}>-\infty$, which can be done via Birkhoff ergodic theorem as in Theorem~\ref{expsepL}, and finishes the proof.
\end{proof}
\subsection{Case 3} In this subsection we briefly explain what happens for the separable Banach space of absolutely continuous functions $AC=AC([-1,0],\RR^N)$  with the Sobolev type norm
\[
\normAC{u}=\normC{u}+\norm{u'}_1=\sup_{s\in[-1,0]} \norm{u(s)}+\int_{-1}^0 \norm{u'(s)}\,ds \quad \text{ for any } u\in AC\,.
\]
We maintain conditions~\ref{S1}, \ref{S2+b} and  \ref{S3}--\ref{S5} of Case 2. The problem now is that this norm is not monotone and the cone is not normal, so that we cannot directly apply the results of  the previous sections. As before,  under assumptions \ref{S1}, \ref{S2+b},  for each $u\in AC$ the initial value problem~\eqref{eq:IVP-C} has a unique Carathéodory type solution  which will be denoted by $z(\cdot,\w,u)$. Moreover, the cocycle relation~\eqref{cocycle} is also satisfied and we can define the linear operator
$ U^{(AC)}_{\w}(t)\colon  AC \longrightarrow  AC$, $u  \mapsto  z_t(\w,u)$. \par\smallskip
First we show that a \mlsps\ is also obtained. Then, although the cone is not normal, we will maintain Definitions~\ref{generalized-floquet-space} and~\ref{def:exponential-separation} without this assumption to prove, by using some of the results in Case 2, the existence of a family of generalized principal Floquet subspaces and a generalized exponential separation of type \textup{II} in this case.  The conclusion of this subsection is that these measurable notions can appear and be used in some natural phase spaces even with weaker properties on the positive cone.
\begin{prop}\label{prop:delay-semiflow-AC}
Under \textup{\ref{S1}} and \textup{\ref{S2+b}}, $\bigl((U^{(AC)}_{\w}(t))_{\w \in \Omega, t \in \mathbb{R}^{+}}, (\theta_{t})_{t \in \R}\bigr)$ is a measurable linear skew\nobreakdash-\hspace{0pt}product semiflow on $AC$ covering $\theta$.
\end{prop}
\begin{proof} From~\cite[Theorem 4.6]{LongoNovoObayaDCDS} we deduce that  $U^{(AC)}_\w(t)\in  \mathcal{L}(AC)$ for each $t\geq 0$ and $\w\in\Omega$. In order to check the $(\mathfrak{B}(\R^{+}) \otimes \mathfrak{F} \otimes \mathfrak{B}(AC), \mathfrak{B}(AC))$\nobreakdash-\hspace{0pt}measurability of the mapping $[\, \R^{+} \times \Omega \times AC \ni (t,\w,u) \mapsto U^{(AC)}_{\w}(t)\,u \in AC \,]$, as shown in~\cite[Lemma 4.51, pp. 153]{AliB}, it is enough to check  that it is a Carathéodory function, i.e. for each $\w\in\Omega$, the map $\R^+\times AC\to AC$, $(t,u) \mapsto  U_\w^{(AC)}(t)\,u=z_t(\w,u)$ is continuous, which follows from~\cite[Theorem 4.6]{LongoNovoObayaDCDS}, and for each $t\in \R^+$ and $u\in AC$
\begin{equation}\label{measurablemapAC}
	[\, \Omega\ni \w  \mapsto  U_\w^{(AC)}(t)\,u=z_t(\w,u)\in AC\,]  \text{ is } (\mathfrak{F}, \mathfrak{B}(AC))\text{\nobreakdash-\hspace{0pt}measurable}.
\end{equation}
Moreover, from the cocycle property~\eqref{eq-cocycle} it is enough to check this property for each $t\in (0,1]$. First  notice that  $AC\simeq \R^N \times L_1([-1,0],\R^N)$ and then its dual space $(AC)^*\simeq \R^N\times L_\infty([-1,0],\R^N)$. Since $AC$ is separable, from Pettis' Theorem (see~Hille and Phillips~\cite[Theorem 3.5.3 and Corollary 2 on pp. 72--73]{HiPhi}) the weak and strong measurability notions are equivalent. Thus, to prove~\eqref{measurablemapAC} it is enough to check that, for  each  $(v_1, v_2)\in \R^N\times L_\infty([-1,0],\R^N)$, $u\in AC$   and $t\in(0,1]$,  the mapping
\begin{equation*}
\bigl[\Omega\ni\w\mapsto \big(\langle v_1,z(t,\w,u)\rangle + \langle v_2,  (U_\w^{(AC)}(t)\,u)'\rangle\big) \in \R\bigr] \text{ is } (\mathfrak{F},
\mathfrak{B}(\R))\text{-measurable}\,,
\end{equation*}
which follows from
\begin{align*}
 \langle  v_2,(U_\w^{(AC)}(t)\,u)'\rangle &=\int_{-1}^0 (v_2(\tau))^t z'(t+\tau,\w,u)\,d\tau \\
  &=\int_{-1}^{-t}   (v_2(\tau))^t\,u'(t+\tau)\, d\tau
    +\int_{0}^t  (v_2(\tau-t))^t \, z'(\tau,\w,u)\,d\tau \\
  &= \int_{-1}^{-t} (v_2(\tau))^t \, u'(t+\tau)\, d\tau +\int_0^t (v_2(\tau-t))^t  A(\theta_{\tau}\w)\, z(\tau,\w,u)\, d\tau \\
 &\quad + \int_0^t (v_2(\tau-t))^t \, B(\theta_{\tau}\w)\,z(\tau-1,\w,u) \,d\tau
\end{align*}
and similar arguments to those of~\cite[Lemmma 4.13]{MiNoOb2}.
\end{proof}
In order to relate both semiflows  $\Phi^{(AC)}$ and $\Phi^{(C)}$ (case 2), we  define for each $t\geq 1$ the linear map
\begin{equation}\label{5.defiskpdCAC}
	\begin{array}{lccc}
		U^{(C,AC)}_\w(t)\colon & C &\longrightarrow & AC\\[.1cm]
		& u & \mapsto & z_t(\w,u)\,.
	\end{array}
\end{equation}
\begin{prop} \label{L(C,AC)}
	Under \textup{\ref{S1}} and \textup{\ref{S2+b}},  for any   $\w\in\Omega$,  $U^{(C,AC)}_\w(t)\in \mathcal{L}(C, AC)$  for each $t\geq 1$  and  it is a compact operator for each $t\geq 2$.
\end{prop}
\begin{proof}
First we check the result for $t=1$ and each $\w\in\Omega$.  Take $u\in AC $; from the definition of the norm,  Proposition~\ref{prop:delay-semiflow-C}  and  the delay differential equation~\eqref{main-delay-eq} we deduce that
\begin{align}  \label{normU)L,C)(1)} \nonumber
	\normAC{U^{(C,AC)}_\w(1)\,u} &=  \normC{z_1(\w,u)}+\int_{-1}^0 \norm{z'(1+s,\w,u)} \, ds \\ \nonumber
	& \leq \normC{z_1(\w,u)}\left(1 + \int_0^1 \!\!a(\theta_s\w) \,ds\right) + \widetilde d(\w) \,\normC{u} \\  \nonumber
	& \leq \left( c(\w)(1+\widetilde{d}(\w))\left(1+ \int_0^1 \!\!a(\theta_s\w) \,ds\right)+ \widetilde{d}(\w)\right)\normC{u}\\
    &:=e(\w) \,\normC{u}
\end{align}
which shows that $U^{(C,AC)}_\w(1)\in \mathcal{L}(C, AC)$. Next, from the cocycle property~\eqref{cocycle} we deduce that  $U_\w^{(C,AC)}(t)= U_{\theta_1\w}^{(C)}(t-1)\circ U_\w^{(C,AC)}(1)\in \mathcal{L}(C, AC)$ for  each $t\geq 1$ and $\w\in\Omega$.
Moreover,  from Lemma~\ref{compactUCS1+S2+b}  we know that $U_\w^{(C)}(t)$ is compact for each $\w\in\Omega$  and $t\geq 1$,  so that the above composition proves  the  compactness of $U_\w^{(C,AC)}(t)$   for each $\w\in\Omega$  and $t\geq 2$, as claimed.
\end{proof}
\begin{coro}
	Under assumptions \textup{\ref{S1}} and \textup{\ref{S2+b}},  $U^{(AC)}_\omega(t)$ is a  compact operator for any $t \ge 2$ and $\omega \in \Omega$.
\end{coro}
\begin{proof}
It is a consequence of  the continuity of  the map $i\colon AC\to C$, $u\mapsto u$ and the previous result.
\end{proof}
As a consequence, from~\cite[Theorem 3.4]{MiNoOb2} we deduce that $\Phi^{(AC)}$ also admits an Oseledets decomposition and the next theorem shows the coincidence of the Lyapunov exponents for both semiflows $\Phi^{(AC)}$ and $\Phi^{(C)}$ (case 2).  We refer the reader to~\cite{MiNoOb2} for all definitions and also the corresponding results for $\Phi^{(C)}$ and $\Phi^{(L)}$ that were proved in that paper.  First we need the following lemma.
\begin{lema}\label{lneL1}
 Assume \textup{\ref{S1}}, \textup{\ref{S2+b}} and consider the function $e$ defined in~\eqref{normU)L,C)(1)}. Then for $\PP$-a.e.\ $\w\in\Omega$
\[ \limsup_{t\to\infty} \frac{\ln e(\theta_t \w)}{t}\leq 0\,.
\]
\end{lema}
\begin{proof}
As in \cite[Lemma 5.1]{MiNoOb2} it is enough to check that $\lnplus{e} \in L^1\OFP$. From~\cite[Lemma 5.6]{MiNoOb1} and inequality~\eqref{ineqc},   we deduce that there is an integer number $n_0$ such that
\[
\lnplus e(\w) \leq 2 \int_0^1 a(\theta_s\w)\,ds + 2\,\lnplus\widetilde{d}(\w) + \ln n_0
\]
and the claim follows from \ref{S2+b}, Fubini's theorem and the invariance of $\PP$.
\end{proof}
\begin{teor}\label{equalLyapunov}
Under assumptions \textup{\ref{S1}} and \textup{\ref{S2+b}},
the sets of Lyapunov exponents for $\Phi^{(AC)}$  and  $\Phi^{(C)}$ coincides.
\end{teor}
\begin{proof}
Let $\lambda^{(C)}$ be a  Lyapunov exponent for $\Phi^{(C)}$ and nonzero $u\in C$. First notice that for each $t\geq 0$  we have $U_\w^{(C)}(t+1)\,u= U_{\theta_t\w}^{(C)}(1)\big(U_\w^{(C)}(t)\,u\big)\in AC$. Thus, from~\eqref{normU)L,C)(1)} we deduce that
\[ \normC{U_\w^{(C)}(t+1)\,u}\leq \normAC{U_\w^{(C,AC)}(t+1)\,u}\leq e(\theta_t\w) \, \normC{U^{(C)}_{\omega}(t) \, u}\,.\]  Therefore, from the definition of $\lambda^{(C)}$ and Lemma~\ref{lneL1} we deduce that
\begin{align*}
\lambda^{(C)}&=\lim_{t\to\infty} \frac{\ln \normC{U_\w^{(C)}(t+1)\,u}}{t+1} \leq \lim_{t\to\infty} \frac{\ln \normAC{U_\w^{(C,AC)}(t+1)\,u}}{t+1} \\
&\leq \limsup_{t \to \infty}\frac{\ln e(\theta_{t}\omega)}{t} + \lim_{t \to \infty} \frac{\ln{\normC{U^{(C)}_{\omega}(t) \,u}}}{t}
 \leq \lambda^{(C)}\,,
\end{align*}
That is,
\[
\lambda^{(C)}= \lim_{t\to\infty} \frac{\ln \normAC{U_\w^{(C,AC)}(t+1)\,u}}{t+1}= \lim_{t\to\infty} \frac{\ln \normAC{U_{\theta_1\w}^{(AC)}(t)\, U_\w^{C}(1)\,u}}{t}\,
\]
is a Lyapunov exponent for $\Phi^{(AC)}$  and nonzero $ U_\w^{C}(1)\,u\in AC$.
When $u\in AC$, the above inequalities show that each Lyapunov exponent for $\Phi^{(AC)}$  and nonzero $u\in AC$  is a Lyapunov exponent for $\Phi^{(C)}$, which finishes the proof.
\end{proof}
\begin{lema}\label{omegaxC->AC} Let  $U_\w^{(C,AC)}(1)$ be the compact linear operator defined in~\eqref{5.defiskpdCAC}. Then  the mapping
	\[
	[\, \Omega\times C\ni (\w,u) \mapsto  U_\w^{(C,AC)}(1)\,u \in AC \,] \text{ is }(\mathfrak{F}\otimes\mathfrak{B}(C), \mathfrak{B}(AC))\text{\nobreakdash-\hspace{0pt}measurable}\,.
	\]
\end{lema}
\begin{proof}
As in Lemma~\ref{omegaxL->C}, it is enough to check that  it is a Carathéodory function,  i.e. for each fixed $\w\in\Omega$, the map $C\to AC$, $u\mapsto  U_\w^{(C,AC)}(1)\,u$ is continuous, which follows from Proposition~\ref{L(C,AC)}, and for each fixed $u\in C$, the map
\begin{equation}\label{measurablemapAC2}
	[\, \Omega\ni \w  \mapsto  U_\w^{(C,AC)}(1)\,u=z_1(\w,u)\in AC\,]  \text{ is } (\mathfrak{F}, \mathfrak{B}(AC))\text{\nobreakdash-\hspace{0pt}measurable}.
\end{equation}
In order to  prove this, we  take a sequence of  functions $u_n\in AC$  converging  to $u\in C$  as $n \uparrow\infty$.  Thus,  the measurability of the maps 	$[\, \Omega\ni \w  \mapsto  z_1(\w,u_n)\in AC\,]$    for each $n\in\N$,  shown in~\eqref{measurablemapAC},  the convergence of them to the map~\eqref{measurablemapAC2} as  $n \uparrow\infty$
and~\cite[Corollary 4.29]{AliB} show the measurability, which finishes the proof.
\end{proof}
\begin{teor}
	Under assumptions \textup{\ref{S1}}, \textup{\ref{S2+b}}, \textup{\ref{S3}} and \textup{\ref{S4}} there is a family of generalized principal Floquet subspaces for the \mlsps\ $\Phi^{(AC)}=\bigl((U^{(AC)}_{\w}(t))_{\w \in \Omega, t \in \RR^+}, (\theta_{t})_{t \in \R}\bigr)$ generated by~\eqref{main-delay-eq}, with generalized principal Lyapunov exponent  $\widetilde{\lambda}^{(AC)}_1=\widetilde{\lambda}^{(C)}_1>-\infty$.
\end{teor}
\begin{proof}
Theorem~\ref{generalized-case2} states the existence of a family $\{E_1^{(C)}(\w)\}_{\w \in \widetilde\Omega_1}$ of generalized principal Floquet subspaces for $\Phi^{(C)}$ with generalized principal Lyapunov exponent  $\widetilde{\lambda}^{(C)}_1$.  For $\w \in \widetilde{\Omega}_1$, put $E_1^{(C)}(\w) := \spanned\{w^{(C)}(\w)\}$ with $\normC{w^{(C)}(\w)} = 1$. First notice that from~\eqref{5.defiskpdCAC}
	\[
	U_{\theta_{-1}w}^{(C,AC)}(1)\,w^{(C)}(\theta_{-1}\w)
	=(z_1(\theta_{-1}\w,w^{(C)}(\theta_{-1}\w)) \quad \text{belongs to } AC
	\]
	and is proportional to $w^{(C)}(\w)$ because $U_{\theta_{-1}w}^{(C)}(1)E_1^{(C)}(\theta_{-1} w)=E_1^{(C)}(\w)$.
Next, we define    $w^{(AC)}(\w) := U_{\theta_{-1}\w}^{(AC)}(1) \, w^{(C)}(\theta_{-1}\w) / \normAC{U_{\theta_{-1}\w}^{(AC)}(1) \, w^{(C)}(\theta_{-1}\w)} $,  unitary and  proportional to $w^{(C)}(\w)$.
	Then $E_1^{(AC)}(\w):=\spanned\{w^{(AC)}(\w)\}=E_1^{(C)}(\w)$  is a one-dimensional subspace of $AC$ and we claim that
	$\{E_1^{(AC)}(\w)\}_{\w \in \widetilde\Omega_1}$  is a family of generalized principal Floquet subspaces for  $\Phi^{(AC)}$, i.e., conditions (i)--(iv) of  Definition~\ref{generalized-floquet-space} are satisfied. First of all, from Lemma~\ref{omegaxC->AC} it is easy to check that $w^{(AC)}\colon\widetilde\Omega_1\to AC^+\setminus \{0\}$ is  $(\mathfrak{F}, \mathfrak{B}(AC))$\nobreakdash-\hspace{0pt}measurable and  thus, (i) holds. Condition~(ii) follows from $U_\w^{(AC)}(t)\,E_1^{(AC)}(\w)=U_\w^{C}(t)\,E_1^{(C)}(\w)=E_1^{(C)}(\theta_t\w)$. From Theorem~\ref{equalLyapunov} we deduce conditions~(iii)--(iv) and the coincidence of the generalized Lyapunov exponents   because both $w^{(C)}(\w)$ and $w^{(AC)}(\w)$ belong to $AC$.
\end{proof}
Finally, once we have checked that $\widetilde\lambda_1^{(AC)}>-\infty$ and $U^{(AC)}_\omega(2)$ is a compact operator, from Theorem~\ref{thm-Oseledets} and the previous theorem an exponential separation of type \textup{II} for the \mlsps\ $\Phi^{(AC)}$~is~obtained.
\begin{teor}
   Assume \textup{\ref{S1}}, \textup{\ref{S2+b}}, \textup{\ref{S3}} and  \textup{\ref{S5}}. Then the \mlsps\ $\Phi^{(AC)}=\bigl((U^{(AC)}_{\w}(t))_{\w \in \Omega, t \in \R^+}, (\theta_{t})_{t \in \R}\bigr)$ generated by system~\eqref{main-delay-eq} admits a generalized exponential separation of type \textup{II} with~{$\widetilde{\lambda}_1^{(AC)}> - \infty$.}
\end{teor}
\begin{proof} From Theorem~\ref{thm-Oseledets}  we obtain a measurable decomposition of $AC$. In addition, from the previous theorem the first term is the one-dimensional subspace $E_1^{(AC)}(\w)$ with
$\lambda_{\mathrm{top}}=\widetilde{\lambda}^{(AC)}_1>-\infty$.  Then, we denote $AC=E_1^{(AC)}(\w)\oplus F_1^{(AC)}(\w)$ where $F_1^{(AC)}(\w)= AC\cap F_1^{(C)}(\w)$  for $\PP$-a.e.\ $\w\in\Omega$. Properties (i) and (ii) of Definition~\ref{def:exponential-separation} are immediate  and finally (iii) is deduced from Theorem~\ref{expIIcase2} following the arguments of Theorem~\ref{expsepC}.
\end{proof}

\end{document}